\documentclass[oneside,11pt]{article}

\usepackage[utf8]{inputenc}
\usepackage[T1]{fontenc}
\usepackage{amsthm}
\usepackage{amsmath}
\usepackage{amsfonts}
\usepackage{amssymb}
\usepackage{bbm}
\usepackage{graphicx} 
\usepackage{float}
\usepackage{booktabs}
\usepackage{multirow}
\usepackage{rotating}
\usepackage{array}
\usepackage{xcolor}
\usepackage[ruled]{algorithm2e}
\usepackage[labelformat=simple]{subcaption}

\DeclareCaptionLabelFormat{subcaptionlabel}{\normalfont(\textbf{#2}\normalfont)}
\newcolumntype{R}[1]{>{\raggedleft\let\newline\\\arraybackslash\hspace{0pt}}m{#1}}
\newcolumntype{L}[1]{>{\raggedright\let\newline\\\arraybackslash\hspace{0pt}}m{#1}}
\usepackage{breakcites}
\usepackage{nowidow}
\usepackage{enumitem}
\usepackage[top=1.5cm,bottom=2cm,right=2.5cm,left=2.5cm]{geometry}
\usepackage{titling}
\usepackage{natbib}
\usepackage{ifplatform}
\usepackage{textcomp}
\ifwindows
\fi
\allowdisplaybreaks

\newcommand{\R}{\mathbb{R}} 
\newcommand{\cF}{\mathcal{F}} 
\newcommand{\cD}{\mathcal{D}} 
\newcommand{\cC}{\mathcal{C}} 
\newcommand{\cB}{\mathcal{B}} 
\newcommand{\cR}{\mathcal{R}} 
\newcommand{\cN}{\mathcal{N}} 
\newcommand{\rmd}{\mathrm{d}} 
\newcommand{\rmb}{\mathsf{b}} 
\newcommand{\rK}{\mathsf{K}} 
\newcommand{\rPhi}{\mathsf{\Phi}} 
\newcommand{\rphi}{\mathsf{\boldsymbol{\phi}}} 
\newcommand{\InfGen}{\mathbb{L}} 
\newcommand{\rmL}{\mathsf{L}} 
\newcommand{\Ind}{\mathbbm{1}} 
\newcommand{\lp}{\left(} 
\newcommand{\rp}{\right)} 
\newcommand{\lb}{\left[} 
\newcommand{\rb}{\right]} 
\newcommand{\lc}{\left\{} 
\newcommand{\rc}{\right\}} 
\newcommand{\blp}{\big(}
\newcommand{\brp}{\big)}

\newcommand{\lrav}[1]{\left|#1\right|} 
\newcommand{\wt}[1]{\widetilde{#1}}

\newcommand{\ol}[1]{\overline{#1}}
\newcommand{\ul}[1]{\underline{#1}}
\newcommand{\lrp}[1]{\lp#1\rp} 
\newcommand{\lrb}[1]{\lb#1\rb} 
\newcommand{\lrc}[1]{\lc#1\rc} 
\newcommand{\blrp}[1]{\blp#1\brp} 
\newcommand{\Pro}[2]{\mathbb{P}_{#2}\lp #1\rp} 
\renewcommand{\Pr}{\mathbb{P}} 
\newcommand{\PPr}{\mathsf{P}} 
\newcommand{\Esp}[1]{\mathbb{E}\lb #1\rb} 
\newcommand{\Espbigg}[1]{\mathbb{E}\bigg[\!#1\bigg]} 
\newcommand{\Es}[2]{\mathbb{E}_{#2}\lb #1\rb} 
\newcommand{\Esbigg}[2]{\mathbb{E}_{#2}\bigg[\!#1\bigg]} 
\newcommand{\E}{\mathbb{E}} 
\newcommand{\EEsp}[1]{\mathsf{E}\lb #1\rb} 
\newcommand{\EEs}[2]{\mathsf{E}_{#2}\lb #1\rb} 
\newcommand{\EE}{\mathsf{E}} 
\newcommand{\VVs}[2]{\mathrm{V}\mathrm{ar}_{#2}\lb #1\rb} 
\newcommand{\Cov}[2]{\mathbb{C}\mathrm{ov}\lb #1,#2\rb} 
\newcommand{\CCov}[2]{\mathsf{C}\mathrm{ov}\lb #1,#2\rb} 

\newtheorem{theorem}{Theorem}

\newtheorem{remark}{Remark}
\newtheorem{proposition}{Proposition}
\newtheorem{lemma}{Lemma}

\newtheorem{example}{Example} 

\usepackage[pdftex,final,bookmarksnumbered,bookmarksopen=false,breaklinks,colorlinks]{hyperref}
\hypersetup{
	final,
	bookmarksnumbered=true,
	bookmarksopen=true,
	bookmarksopenlevel=0,
	unicode=false,
	pdftoolbar=true,
	pdfmenubar=true,
	pdffitwindow=false,
	pdftitle={Optimal stopping of Gauss-Markov bridges},
	pdfdisplaydoctitle=true,
	pdftoolbar=true,
	pdfmenubar=true,
	pdflang={English},
	pdfauthor={Abel Azze, Bernardo D'Auria, Eduardo Garcia-Portugues},
	pdfsubject={arXiv paper},
	pdfcreator={Abel Azze},
	pdfproducer={Abel Azze},
	pdfkeywords={Brownian bridge}{Gauss–Markov bridge}{Optimal stopping}{Ornstein–Uhlenbeck bridge}{Time-inhomogeneity},
	pdfnewwindow=true,
	breaklinks=true,
	hidelinks,       
	linkcolor=black,
	citecolor=black,
	filecolor=magenta,
	urlcolor=cyan
}

\newif\ifmain
\maintrue
\newif\ifsupplement
\supplementtrue

\newif\iffigstabs
\figstabstrue

\begin{document}

\ifmain

\title{Optimal stopping of Gauss--Markov bridges}
\setlength{\droptitle}{-1cm}
\predate{}%
\postdate{}%
\date{}

\author{Abel Azze$^{1,4}$, Bernardo D'Auria$^{2}$, and Eduardo Garc\'ia-Portugu\'es$^{3}$}
\footnotetext[1]{Department of Quantitative Methods, CUNEF Universidad (Spain).}
\footnotetext[2]{Department of Mathematics ``Tullio Levi Civita'', University of Padova (Italy).}
\footnotetext[3]{Department of Statistics, Universidad Carlos III de Madrid (Spain).}
\footnotetext[4]{Corresponding author. e-mail: \href{mailto:abel.guada@cunef.edu}{abel.guada@cunef.edu}.}
\maketitle

\begin{abstract}
	We solve the non-discounted, finite-horizon optimal stopping problem of a Gauss--Markov bridge by using a time-space transformation approach. The associated optimal stopping boundary is proved to be Lipschitz continuous on any closed interval that excludes the horizon, and it is characterized by the unique solution of an integral equation. A Picard iteration algorithm is discussed and implemented to exemplify the numerical computation and geometry of the optimal stopping boundary for some illustrative cases.
\end{abstract}
\begin{flushleft}
\end{flushleft}

\section{Introduction}\label{sec:intro}

The problem of optimally stopping a Markov process to attain a maximum mean reward dates back to Wald's sequential analysis \cite{Wald_1947_sequential} and is consolidated in the work of \cite{Dynkin_1963_optimum}. Ever since, it has received increasing attention from numerous theoretical and practical perspectives, as comprehensively compiled in the book of \cite{Peskir_2006_optimal}. However, Optimal Stopping Problems (OSPs) are mathematically complex objects, which makes it difficult to obtain sound results in general settings, and typically lead to requiring smoothness conditions and simplifying assumptions for their solution. One of the most popular simplifying assumptions is the time-homogeneity of the underlying Markovian process.

Time-inhomogeneous diffusions can be cast back to time-homogeneity (see, e.g., \cite{Taylor_1968_optimal}, \cite{Dochviri_1995_optimal}, \cite{Shiryaev_2008_optimal}) at the cost of increasing the dimension of the OSP, which results in an increased complexity, hampering subsequent derivations or limiting studies to tackle specific, simplified time dependencies. Take, as examples, the works of \cite{Krylov_1980_controlled}, \cite{Oshima_2006_optimal}, and \cite{Yang_2014_refined}, which proved different types of continuities and characterizations of the value function; or those of \cite{Friedman_1975_stopping} and \cite{Jacka_1992_finite-horizon}, which shed light on the shape of the stopping set; and \cite{Friedman_1975_parabolic} and \cite{Peskir_2019_continuity}, who studied the smoothness of the associated free boundary. To mitigate the burden of time-inhomogeneity, many of these works ask for the process' coefficients to be Lipschitz continuous or at least bounded. This usual assumption excludes important classes of time-dependent processes, such as diffusion bridges, whose drifts explode as time approaches a terminal point.

In a broad and rough sense, bridge processes, or bridges for short, are stochastic processes ``anchored'' to deterministic values at some initial and terminal time points. Formal definitions and potential applications of different classes of bridges have been extensively studied. Bessel and Lévy bridges are respectively described by \cite{Pitman_1982_decomposition} and \cite{Salminen_1984_Brownian}, and by \cite{Hoyle_2011_Levy} and \cite{Erickson_2020_probability}. A canonical reference for Gaussian bridges can be found in the work of \cite{Gasbarra_2007_Gaussian}, while Markov bridges are addressed in great generality by \cite{Fitzsimmons_1993_Markovian}, \cite{Chaumont_2011_Markovian}, and \cite{Cetin_2016_Markov}.

In finance, diffusion bridges are appealing models from the perspective of a trader who wants to incorporate his beliefs about future events, like in trading perishable commodities, modeling the presence of arbitrage, incorporating algorithms' forecasts and experts' predictions, or trading mispriced assets that could rapidly return to their fair price. Works that consider models based on a Brownian Bridge (BB) to address these and other insider trading situations include \cite{Kyle_1985_continuous}, \cite{Brennan_1990_arbitrage}, \cite{Back_1992_insider}, \cite{Liu_2004_losing}, \cite{Campi_2007_insider}, \cite{Campi_2011_dynamic}, \cite{Campi_2013_equilibrium}, \cite{Cetin_2013_point}, \cite{Sottinen_2014_generalized}, \cite{Cartea_2016_algorithmic}, \cite{Angoshtari_2019_optimal}, and \cite{Chen_2021_constrained}. The early work of \cite{Boyce_1970_stopping} had already suggested the use of a BB to model the perspective of an investor who wants to optimally sell a bond. Recently, \cite{DAuria_2020_discounted} applied a BB to optimally exercise an American option in the presence of the so-called stock-pinning effect (see \cite{Nelken_2001_effect}, \cite{Ni_2005_stock}, \cite{Golez_2012_pinning}, and \cite{Ni_2021_does}), obtaining competitive empirical results when compared to the classic Black--Scholes model.
Taking distance from the BB model, \cite{Hilliard_2015_pricing} used an Ornstein--Uhlenbeck Bridge (OUB) to model the effect of short-lived arbitrage opportunities in pricing an American option, recurring to a binomial-tree numerical method instead of deriving analytical results.

Non-financial applications of BBs include their usual adoption to model animal movement (see \cite{Horne_2007_analyzing, Venek_2016_evaluating}, \cite{Kranstauber_2019_modelling}, and \cite{Krumm_2021_Brownian}), and their construction as a limit case of sequentially drawing elements without replacement from a large population (see \cite{Rosen_1965_limit}). The latter connection makes BBs good asymptotic models for classical statistical problems, like variations of the urn's problem (see \cite{Ekstrom_2009_optimal}, \cite{Andersson_2012_card}, and \cite{Chen_2015_optimal}).

Whenever the goal is to optimize the time to take an action, all the previous situations in which a BB, an OUB, or diffusion bridges have applications can be intertwined with optimal stopping theory. However, within the time-inhomogeneous realm, diffusion bridges are particularly challenging to treat with classical optimal-stopping tools, as they feature explosive drifts. It comes as no surprise, hence, that the literature addressing this topic is scarce when compared with its non-bridge counterpart. The first incursion into OSPs with diffusion bridges is by Shepp's work \cite{Shepp_1969_explicit}, who solved the OSP of a BB by linking it to that of a simpler Brownian Motion (BM) representation. After Shepp's result, the more recent studies of OSPs with diffusion bridges still revolve around variations of the BB. \cite{Ekstrom_2009_optimal} and \cite{Ernst_2015_revisiting} revisited Shepp's problem with novel solution methods. \cite{Ekstrom_2009_optimal} and \cite{DeAngelis_2020_optimal} widened the class of gain functions; \cite{DAuria_2020_discounted} considered the (exponentially) discounted version; while \cite{Follmer_1972_optimal}, \cite{Leung_2018_optimal}, \cite{Glover_2020_optimally}, and \cite{Ekstrom_2020_optimal}, introduced randomization in either the terminal time or the pinning point. To the best of our knowledge, the only solution to an OSP with diffusion bridges that steps outside the BB, came recently in \cite{DAuria-2021-optimal}, which extends Shepp's technique to embrace an OUB.

Both the BB and the OUB belong to the class of Gauss--Markov Bridges (GMBs), that is, bridges that simultaneously exhibit the Markovian and Gaussian properties. Due to their enhanced tractability and wide applicability, these processes have been in the spotlight for some decades, especially in recent years. A good compendium of works related to GMBs can be found in \cite{Abrahams_1981_some}, \cite{Buonocore_2013_some}, \cite{Barczy_2013_representation}, \cite{Barczy_2013_sample}, \cite{Barczy_2011_general}, \cite{Chen_2016_stochastic}, and \cite{Hildebrandt_2020_pinned}.

In this paper we solve the finite-horizon OSP of a GMB. In doing so, we generalize not only Shepp's result for the BB case, but also its methodology. Indeed, the same type of transformation that casts a BB into a BM is embedded in a more general change-of-variable method to solve OSPs, which is detailed in \cite[Section 5.2]{Peskir_2006_optimal} and exemplary used in \cite{Pedersen_2002_onnonlinear} for non-linear OSPs. When the GM process is also a bridge, such a representation presents regularities that we show useful to overcome the bridges' explosive drifts. Loosely, the drift's divergence is equated to that of a time-transformed BM and then explained in terms of the laws of iterated logarithms. This trick allows working out the solution of an equivalent infinite-horizon OSP with a time-space transformed BM underneath, and then casting the solution back into original terms. The solution is attained, in a probabilistic fashion, by proving that both the value function and the Optimal Stopping Boundary (OSB) are regular enough to meet the premises of a relaxed Itô's lemma that allows deriving the free-boundary equation. In particular, we prove the Lipschitz continuity of the OSB, which we use to derive the global continuous differentiability of the value function and, consequently, the smooth-fit condition. The free-boundary equation is given in terms of a Volterra-type integral equation with a unique solution. For enriched perspectives and full sight of the reach of GMBs, we provide, besides the BM representation, a third angle from which GMBs can be seen: as time-inhomogeneous OUBs. Hence, our work also extends the work of \cite{DAuria-2021-optimal} for a time-independent OUB. This OUB representation is arguably more appealing to numerically explore the OSB's shape, which is done by using a Picard iteration algorithm that solves the free-boundary equation. The OSB exhibits a trade-off between two pulling forces, the one towards the mean-reverting level of the OUB representation, and that which anchors the process at the horizon. The numerical results also reveal that the OSB is not monotonic in general, making this paper one of the few results in the optimal-stopping literature that characterizes non-monotonic OSBs in a general framework.

The rest of this paper is organized as follows. Section \ref{sec:GMB} establishes four equivalent definitions of GMBs, including the time-space transformed BM representation. Section \ref{sec:formulation} introduces the finite-horizon OSP of a GMB and proves its equivalence to that of an infinite-horizon, time-dependent gain function, and a BM underneath. The auxiliary OSP is then treated in Section \ref{sec:solution} as a standalone problem. This section also accounts for the main technical work of the paper, where classical and new techniques of optimal stopping theory are combined to obtain the solution of the OSP. This solution is then translated back into original terms in Section \ref{sec:solutionoriginal}, where the free-boundary equation is provided. Section \ref{sec:numerical_results} discusses the practical aspects of numerically solving the free-boundary equation, and shows computer drawings of the OSB. Final remarks are given in Section~\ref{sec:conclusions}.

\section{Gauss--Markov bridges}\label{sec:GMB}

Both Gaussian and Markovian processes exhibit features that are appealing from a theoretical, computational, and applicable viewpoint. Gauss--Markov (GM) processes, that is, processes that are Gaussian and Markovian at the same time, merge the advantages of these two classes. They inherit the convenient Markovian lack of memory and the Gaussian processes' property of being characterized by their mean and covariance functions. Additionally, the Markovianity of Gaussian processes is equivalent for their covariances to admit a certain ``factorization''. The following lemma collects such a useful characterization, whose proof follows from the lemma on page 863 from \cite{Borisov_1983_criterion}, and Theorem 1 and Remarks 1--2 in \cite{Mehr_1965_certain}.

Here and thereafter, when we mention a non-degenerated GM process in an interval, we mean that its marginal distributions are non-degenerated in the same interval. In addition, we always consider the GM processes defined in their natural filtrations.

\begin{lemma}[Characterization of non-degenerated GM processes]\label{lm:GM_def}\ \\
    A function $R:[0, T]^2\rightarrow\R$ such that $R(t_1, t_2) \neq 0$ for all $t_1, t_2\in(0, T)$ is the covariance function of a non-degenerated GM process in $(0, T)$ if and only if there exist functions $r_1, r_2:[0, T]\rightarrow\R$, that are unique up to a multiplicative constant, such that
    \begin{enumerate}[label=(\textit{\roman{*}}), ref=(\textit{\roman{*}})]
        \item $R(t_1, t_2) = r_1(t_1\wedge t_2)r_2(t_1\vee t_2)$; \label{eq:cov_fac}
        \item $r_1(t) \neq 0$ and $r_2(t) \neq 0$ for all $t\in(0, T)$; \label{eq:cov_non-zero}
        \item $r_1/r_2$ is positive and strictly increasing on $(0, T)$. \label{eq:cov_ratio}
    \end{enumerate}
    Moreover, $r_1$ and $r_2$ take the form
    \begin{align}\label{eq:explicit_cov_fac}
        r_1(t) = 
        \begin{cases}
        R(t, t'), & t\leq t', \\
        R(t, t)R(t', t')/R(t', t), & t > t',
        \end{cases} \quad 
        r_2(t) = 
        \begin{cases}
        R(t, t)/R(t, t'), & t\leq t', \\
        R(t', t)/R(t', t'), & t > t',
        \end{cases}
    \end{align}
    for some $t'\in(0, T)$. Changing $t'$ is equivalent to scaling $r_1$ and $r_2$ by a constant factor.
\end{lemma}

We say that the functions $r_1$ and $r_2$ in Lemma \ref{lm:GM_def} are a factorization of the covariance function $R$. The lemma provides a simple technique to construct GM processes with ad hoc covariance functions that are not necessarily time-homogeneous. This is particularly useful given the complexity of proving the positive-definiteness of an arbitrary function to check its validity as a covariance function. GM processes also admit a simple representation by means of time-space transformed BMs (see, e.g., \cite{Mehr_1965_certain}), which results in higher tractability. Moreover, viewed through the lens of diffusions, GM processes account for space-linear drifts and space-independent volatilities, both coefficients being time-dependent (see, e.g., \cite{Buonocore_2013_some}).

We call Gauss--Markov Bridge (GMB) a process that results after ``conditioning'' (see, e.g., \cite{Gasbarra_2007_Gaussian} for a formal definition) a GM process to start and end at some initial and terminal points. It is straightforward to see that the Markovian property is preserved after conditioning. Although not as evidently, the bridge process also inherits the Gaussian property (see, e.g., \cite[Formula A.6]{Ramussen_2006_Gaussian}, or \cite{Buonocore_2013_some}). Hence, the above-mentioned conveniences of GM processes are inherited by GMBs. In particular, the time-space transformed BM representation adopts a specific form that characterizes GMBs and forms the backbone of our main results. The following proposition sheds light on that representation and serves to formally define a GMB as well as to offer different characterizations.

\begin{proposition}[Gauss--Markov bridges]\label{pr:GMB}\ \\
    Let $\smash{X = \{X_u\}_{u\in [0, T]}}$ be a GM process defined on the probability space $(\Omega, \cF, \PPr)$, for some $T > 0$. The following statements are equivalent:
    \begin{enumerate}[label=(\textit{\roman{*}}), ref=\textit{\roman{*}}]
    \item \label{def:GMB_derived} There exists a time-continuous GM process, non-degenerated on $[0, T]$, defined on $(\Omega, \cF, \PPr)$, and denoted by $\smash{\wt{X} = \{\wt{X}_u\}_{u\in[0, T]}}$, whose mean and covariance functions are twice continuously differentiable, and such that
          \begin{align*}
            \mathrm{Law}(X, \PPr) = \mathrm{Law}(\wt{X}, \PPr_{x, T, z}),
          \end{align*}
          with $\PPr_{x, T, z}(\cdot) = \PPr(\cdot\ |\wt{X}_0 = x, \wt{X}_T = z)$ for some $x\in\R$ and $(T, z)\in\R_+\times\R$.
    \item \label{def:GMB_cov} Let $m(t) := \EEsp{X_t}$ and $R(t_1, t_2) := \CCov{X_{t_1}}{X_{t_2}}$, where $\EE$ and $\mathsf{C}\mathrm{ov}$ are the mean and covariance operators related to $\PPr$. Then, $t \mapsto m(t)$ is twice continuously differentiable, and there exist functions $r_1$ and $r_2$ that are unique up to multiplicative constants and such that:
    \begin{enumerate}[label=(\ref{def:GMB_cov}.\textit{\arabic{*}}), ref=\ref{def:GMB_cov}.\textit{\arabic{*}}]
        \item $R(t_1, t_2) = r_1(t_1\wedge t_2)r_2(t_1\vee t_2)$; \label{eq:cov_1}
        \item $r_1(t) \neq 0$ and $r_2(t) \neq 0$ for all $t\in(0, T)$; \label{eq:cov_2}
        \item $r_1/r_2$ is positive and strictly increasing on $(0, T)$; \label{eq:cov_3}
        \item $r_1(0) = r_2(T) = 0$; \label{eq:cov_4}
        \item $r_1$ and $r_2$ are twice continuously differentiable; \label{eq:cov_5}
        \item $r_1(T) \neq 0$ and $r_2(0) \neq 0$. \label{eq:cov_6}
    \end{enumerate}
    \item \label{def:GMB_BM} $X$ admits the representation
    \begin{align}\label{eq:GMB_to_BM}
        \hspace*{-0.5cm}\left\{
        \begin{aligned}
             X_t &= \alpha(t) + \beta_T(t)\lrp{(z - \alpha(T))\gamma_T(t) + \lrp{B_{\gamma_T(t)} + \frac{x - \alpha(0)}{\beta_T(0)}}},\; t\in[0, T), \\
            X_{T} &= z.
        \end{aligned}
        \right.
    \end{align}
    where $\lrc{B_u}_{u\in\R_+}$ is a standard BM, and $\alpha:[0, T]\rightarrow\R$, $\beta_T:[0, T]\rightarrow\R_+$, and $\gamma_T:[0, T)\rightarrow\R_+$ are twice continuously differentiable functions such that: 
    \begin{enumerate}[label=(\ref{def:GMB_BM}.\textit{\arabic{*}}), ref=\ref{def:GMB_BM}.\textit{\arabic{*}}]
        \item $\beta_T(T) = \gamma_T(0) = 0$; \label{pr:GMB_to_BM_1}
        \item $\gamma_T$ is monotonically increasing; \label{pr:GMB_to_BM_2} 
        \item $\lim_{t\rightarrow T} \gamma_T(t) = \infty$ and $\lim_{t\rightarrow T} \beta_T(t)\gamma_T(t) = 1$. \label{pr:GMB_to_BM_3}
    \end{enumerate}
    \item \label{def:GMB_OUB} 
    $X$ is the unique strong solution of the following OUB Stochastic Differential Equation (SDE)
    \begin{align}\label{eq:tdOU_SDE}
        \rmd X_t = \theta(t)(\kappa(t) - X_t)\,\rmd t + \nu(t)\,\rmd B_t,\quad t\in(0, T),
    \end{align}
    with initial condition $X_0 = x$. $\lrc{B_t}_{t\in\R_+}$ is a standard BM, and $\theta:[0, T)\rightarrow\R_+$, $\kappa:[0, T]\rightarrow\R$, and $\nu:[0, T]\rightarrow\R_+$ are continuously differentiable functions such that:
    \begin{enumerate}[label=(\ref{def:GMB_OUB}.\textit{\arabic{*}}), ref=\ref{def:GMB_OUB}.\textit{\arabic{*}}]
        \item $\lim_{t\rightarrow T}\int_0^t\theta(u)\,\rmd u = \infty$; \label{eq:theta_explotion}
        \item $\nu^2(t) = \theta(t)\exp\lrc{-\int_0^t\theta(u)\,\rmd u}$ or, equivalently, $\theta(t) = \nu^2(t)\big/\int_t^T\nu^2(u)\,\rmd u$. \label{eq:theta_sigma}
    \end{enumerate}
    \end{enumerate}
\end{proposition}

\begin{proof}
    \eqref{def:GMB_derived} $\Longrightarrow$ \eqref{def:GMB_cov}. $X$ is a non-degenerated GM process on $(0, T)$, as it arises by conditioning a process with the same qualities to take deterministic values at $t = 0$ and $t = T$. Hence, Lemma \ref{lm:GM_def} guarantees that $R(t_1, t_2) := \CCov{X_{t_1}}{X_{t_2}}$ meets the conditions \eqref{eq:cov_1}--\eqref{eq:cov_3}. Since $X$ degenerates at $t = 0$ and $t = T$, and due to \eqref{eq:cov_1}, condition \eqref{eq:cov_4} holds true. From the twice continuous differentiability (with respect to both variables) of the covariance function of $\wt{X}$, it follows that of $X$ which, alongside \eqref{eq:explicit_cov_fac}, implies \eqref{eq:cov_5}. 
    
    We now prove \eqref{eq:cov_6}. Let $\wt{m}, \wt{r}_1, \wt{r}_2:[0, T]\rightarrow\R$ be the mean and the covariance factorization of $\wt{X}$. Hence (see, e.g., \cite[Formula A.6]{Ramussen_2006_Gaussian} or \cite{Buonocore_2013_some}),
    \begin{align}\label{eq:dGMB_mean}
        m(t) = \wt{m}(t) + (x - \wt{m}(0))\frac{r_2(t)}{r_2(0)} + 
        (z - \wt{m}(T))r_1(t),\quad t\in[0, T),
    \end{align}
    and
    \begin{align}\label{eq:dGMB_cov_fac}
    \left\{
    \begin{aligned}
        r_1(t) &= \frac{\wt{r}_1(t)\wt{r}_2(0) - \wt{r}_1(0)\wt{r}_2(t)}{\wt{r_1}(T)\wt{r}_2(0) - \wt{r}_1(0)\wt{r}_2(T)}, \\
        r_2(t) &= \wt{r}_1(T)\wt{r_2}(t) - \wt{r}_1(t)\wt{r}_2(T).
    \end{aligned}
    \right.
    \end{align}
    From the continuity of $\wt{R}$ and the representation \eqref{eq:explicit_cov_fac}, it follows the continuity of $\wt{r}_1/\wt{r}_2$. Note that $\wt{r}_2$ does not vanish at $t = 0$ and $t = T$ due to the non-degenerated nature of $\wt{X}$ at both boundary points. Hence, we can extend the increasing nature of $\wt{r}_1/\wt{r}_2$, established in \ref{eq:cov_ratio} from Lemma \ref{lm:GM_def}, to $t = 0$ and $t = T$, which implies that $\wt{r}_1(T)\wt{r}_2(0) - \wt{r}_1(0)\wt{r}_2(T) > 0$ and, therefore, \eqref{eq:dGMB_cov_fac} results in $r_1(T)=1$ and $r_2(0) > 0$. This does not mean that $r_1(T)$ and $r_2(0)$ must be positive, as $-r_1$ and $-r_2$ are also a factorization of $R$, but it does imply \eqref{eq:cov_6}.
    
    \eqref{def:GMB_cov} $\Longrightarrow$ \eqref{def:GMB_derived}. Consider the functions
    \begin{align}\label{eq:parent_mean}
        \wt{m}(t) &:= m(t) - (x - m_1)\frac{r_2(t)}{r_2(0)} - (z - m_2)r_1(t), \quad t \in(0, T),
    \end{align}
    with $\wt{m}(0) := m_1$ and $\wt{m}(T) := m_2$ for $m_1, m_2 \in\R$, and
    \begin{align}
        \wt{r}_1(t) := a r_1(t) + b r_2(t); \quad \wt{r}_2(t) := c r_1(t) + d r_2(t), \quad t \in[0, T], \label{eq:parent_cov_fac}
    \end{align}
    for $a, b, c, d > 0$ and such that $ad > bc$. This relation is met, for instance, by setting $a = b = c = 1$ and $d = 2$. We can divide by $r_2(0)$ in \eqref{eq:parent_mean} since \eqref{eq:cov_6} holds true. Let $h(t) := r_1(t)/r_2(t)$ and $\wt{h}(t) := \wt{r}_1(t)/\wt{r}_2(t)$. We get $\wt{h}(t) = (ah(t) + b)/(ch(t) + d)$ from~\eqref{eq:parent_cov_fac}. Hence,
    \begin{align*}
        \wt{h}'(t) > 0 \Longleftrightarrow h'(t)\lrp{ad - bc} > 0.
    \end{align*}
    Condition \eqref{eq:cov_3} along with our choice of $a$, $b$, $c$, and $d$ guarantees that the right-hand side of the equivalence holds. Therefore, $\wt{h}(t)$ is strictly increasing. Since $\wt{h}$ is also positive, $\wt{R}(t_1, t_2) := \wt{r}_1(t_1\wedge t_2)\wt{r}_2(t_1\vee t_2)$ is the covariance function of a non-degenerated GM process, as stated in Lemma \ref{lm:GM_def}. Let $\wt{X} = \{\wt{X}_t\}_{t\in [0, T]}$ be a GM process with mean $\wt{m}(t)$ and covariance $\wt{R}(t_1, t_2)$. From the differentiability of $m$, $r_1$, and $r_2$, alongside \eqref{eq:parent_mean} and \eqref{eq:parent_cov_fac}, it follows that of $\wt{m}$, $\wt{r}_1$, and $\wt{r}_2$ (and $\wt{R}$).
    
    One can check, after some straightforward algebra and in alignment with \eqref{eq:dGMB_mean}--\eqref{eq:dGMB_cov_fac}, that the mean and covariance functions of the GMB derived from conditioning $\wt{X}$ to go from $(0, x)$ to $(T, z)$ coincide with $m$ and $R$.
    
    \eqref{def:GMB_derived} $\Longrightarrow$ \eqref{def:GMB_BM}. Let $\wt{m}(t) := \EE[\wt{X}_t]$ and $\wt{R}(t_1, t_2) := \CCov{\wt{X}_{t_1}}{\wt{X}_{t_2}}$. As a result of conditioning $\wt{X}$ to have initial and terminal points $(0, x)$ and $(T, z)$, $X$ is a GM process with mean $m$ given by \eqref{eq:dGMB_mean} and covariance factorization $r_1$ and $r_2$ given by \eqref{eq:dGMB_cov_fac}. Although not explicitly indicated, recall that $m$ depends on $x$, $T$, and $z$, and $r_1$ and $r_2$ depend on $T$.
    
    Therefore, $X$ admits the representation
    \begin{align}\label{eq:BM_rep}
        X_t = m(t) + r_2(t)B_{h(t)},\quad 0\leq t < T,
    \end{align}
    where $t\mapsto h(t) := r_1(t)/r_2(t)$ is a strictly increasing function such that $h(0) = 0$ and $\lim_{t\rightarrow T}h(t) = \infty$. Since $\lim_{t\rightarrow T}r_2(t)h(t) = r_1(T) = 1$ (see \eqref{eq:dGMB_cov_fac}), the law of the iterated logarithm allows us to continuously extend $X_t$ to $T$ as the $\PPr$-a.s. limit $X_T := \lim_{t\rightarrow T}X_t = z$. Then, representation \eqref{eq:GMB_to_BM} and properties \eqref{pr:GMB_to_BM_1}--\eqref{pr:GMB_to_BM_3} follow after taking $\alpha = \wt{m}$, $\beta_T = r_2$, and $\gamma_T = h$. It also follows that $\alpha$, $\beta_T$, and $\gamma_T$ are twice continuously differentiable, like $\wt{m}$, $\wt{r}_1$, and $\wt{r}_2$ as well.
    
    \eqref{def:GMB_BM} $\Longrightarrow$ \eqref{def:GMB_cov}. Assuming that $X = \{X_t\}_{t\in[0, T]}$ admits representation \eqref{eq:GMB_to_BM} and that properties \eqref{pr:GMB_to_BM_1}--\eqref{pr:GMB_to_BM_3} hold, then $X$ is a GMB with covariance factorization given by $r_1(t) = \beta_T(t_1)\gamma_T(t_1)$ and $r_2(t) = \beta_T(t)$. It readily follows that $r_1$ and $r_2$ meet conditions \eqref{eq:cov_1}--\eqref{eq:cov_6}. It is also trivial to note that $X$ has a twice continuously differentiable mean.
    
    \eqref{def:GMB_derived} $\Longrightarrow$ \eqref{def:GMB_OUB}.
    Let $\EE_{t, x}$ and $\E_{s, y}$ be the mean operators with respect to the probability measures $\PPr_{t, x}$ and $\Pr_{s, y}$, such that $\PPr_{t, x}(\cdot) = \PPr(\cdot | X_t = x)$ and $\Pr_{s, y}(\cdot) = \Pr(\cdot | B_s = y)$, where $\lrc{B_u}_{u\in\R_+}$ is the BM in representation \eqref{eq:BM_rep}. 
    Then,
    \begin{align*}
        \mathrm{Law}\lrp{\lrc{X_u}_{u\in [t, T)}, \PPr_{t, x}} = \mathrm{Law}\lrp{\lrc{m(u) + r_2(u)B_{h(u)}}_{u\in[t, T)}, \Pr_{s, y}},
    \end{align*}
    for $s = h(t)$ and $y = (x - m(t))/r_2(t)$. Hence,
    \begin{align*}
        \EEs{X_{t+\varepsilon} - x}{t, x} &= \EEs{m(t+\varepsilon) + r_2(t+\varepsilon)B_{h(t+\varepsilon)} - x}{s, y} \\
        &= \EEs{m(t+\varepsilon) + \frac{r_2(t+\varepsilon)}{r_2(t)}(x- m(t)) + r_2(t+\varepsilon)B_{h(t+\varepsilon)-h(t)} - x}{s, y}.
    \end{align*}
    Likewise,
    \begin{align*}
        \EEs{(X_{t + u} - x)^2}{t, x} &= \EEsp{\lrp{m(t+\varepsilon) + \frac{r_2(t+\varepsilon)(x - m(t))}{r_2(t)} + r_2(t+\varepsilon)B_{h(t+\varepsilon) - h(t)} - x}^2} \\
        &= \lrp{m(t+\varepsilon) + \frac{r_2(t+\varepsilon)(x - m(t))}{r_2(t)} - x}^2 + r_2^2(t+\varepsilon)(h(t+\varepsilon) - h(t)).
    \end{align*}
    Therefore,
    \begin{align*}
        \lim_{\varepsilon\downarrow 0}\varepsilon^{-1}\EEs{X_{t+\varepsilon} - x}{t, x} &= m'(t) + (x - m(t))r_2'(t)/r_2(t), \\
        \lim_{\varepsilon\downarrow 0}\varepsilon^{-1}\EEs{\lrp{X_{t+\varepsilon} - x}^2}{t, x} &= r_2^2(t)h'(t).
    \end{align*}
    By comparing the drift and volatility terms, $X$ is the unique strong solution (see Example 2.3 by \cite{Cetin_2016_Markov}) of the SDE \eqref{eq:tdOU_SDE} for 
    \begin{align}\label{eq:tdOUB_coeff_mean_cov}
        \left\{
    \begin{aligned}
        \theta(t) &= -r_2'(t)/r_2(t),  \\ 
        \kappa(t) &= m(t) - m'(t)r_2(t)/r_2'(t), \\
        \nu(t) &= r_2(t)\sqrt{h'(t)}. 
    \end{aligned}
        \right.
    \end{align}
    It follows from \eqref{eq:tdOUB_coeff_mean_cov} (or by directly deriving it from \eqref{eq:tdOU_SDE}) that
    \begin{align}
       m(t) &= \varphi(t)\lrp{x + \int_0^t \frac{\kappa(u)\theta(u)}{\varphi(u)}\,\rmd u} \label{eq:tdOU_mean_1} \\
       &= \varphi(t)\lrp{x + \int_0^t \frac{\wt{m}(u)\theta(u) - \wt{m}'(u)}{\varphi(u)}\,\rmd u + (z - \wt{m}(T))\int_0^t \frac{r_1(u)\theta(u) - r_1'(u)}{\varphi(u)}\,\rmd u} \label{eq:tdOU_mean_2}
    \end{align}
    and
    \begin{align}\label{eq:tdOU_cov_fac}
        r_1(t) = \varphi(t)\int_0^t\frac{\nu^2(u)}{\varphi^2(u)}\,\rmd u,\quad r_2(t) = \varphi(t),
    \end{align}
    for $t\in[0, T)$, with $\varphi(t) = \exp\lrc{-\int_0^t \theta(u)\,\rmd u}$. Since $X$ is degenerated at $t = T$, $r_2(T) = 0$, which implies \eqref{eq:theta_explotion}. 
    By comparing \eqref{eq:tdOU_mean_2} with \eqref{eq:dGMB_mean},
    \begin{align*}
        r_1(t) = \varphi(t)\int_0^t \frac{r_1(u)\theta(u) - r_1'(u)}{\varphi(u)}\,\rmd u = 2\varphi(t)\int_0^t \frac{r_1(u)\theta(u)}{\varphi(u)}\,\rmd u - r_1(t),
    \end{align*}
    which, after using \eqref{eq:tdOU_cov_fac}, leads to
    \begin{align*}
        \int_0^t \frac{\nu^2(u)}{\varphi^2(u)}\,\rmd u = \int_0^t \frac{r_1(u)\theta(u)}{\varphi(u)}\,\rmd u.
    \end{align*}
    Differentiating with respect to $t$ both sides of the equation above, and relying again on \eqref{eq:tdOU_cov_fac}, we get
    \begin{align*}
       \frac{\nu^2(t)}{\varphi^2(t)}=\theta(t)\int_0^t\frac{\nu^2(u)}{\varphi^2(u)}\,\rmd u.
    \end{align*}
    The expression above is an ordinary differential equation in $f(t) = \int_0^t\nu^2(u)/\varphi^2(u)\,\rmd u$ whose solution is $f(t) = C_1 + 1/\varphi(t)$ for some constant $C_1$. Hence, $f'(t) = \theta(t)/\varphi(t)$. Therefore, some straightforward algebra leads us to the first equality in \eqref{eq:theta_sigma}, which implies that
    \begin{align*}
       \int_0^t\nu^2(u)\,\rmd u = C_2 + \int_0^t\theta(u)\varphi(u)\,\rmd u = C_2 + 1 - \varphi(t),
    \end{align*}
    for a constant $C_2\in\R$. Since $\lim_{t\rightarrow T}\varphi(t) = 0$, then $C_2 = \int_0^T\nu^2(u)\,\rmd u - 1$. Hence,
    \begin{align*}
        \int_0^t\theta(u)\,\rmd u &= -\ln\lrp{C_2 + 1 - \int_0^t\nu^2(u)\,\rmd u},
    \end{align*}
    from where it follows the second equality in \eqref{eq:theta_sigma} after differentiating.
    
    Finally, from the smoothness of $\wt{m}$, $\wt{r}_1$, and $\wt{r}_2$, which implies that of $m$, $r_1$, and $r_2$, it follows that $\theta$, $\kappa$, and $\nu$ are continuously differentiable.
    
    \eqref{def:GMB_OUB} $\Longrightarrow$ \eqref{def:GMB_cov}. Functions $\theta$, $\kappa$, and $\nu$ are sufficiently regular to prove, by using Itô's lemma, that
    \begin{align*}
        X_t = \varphi(t)\lrp{X_0 + \int_0^t\frac{\kappa(u)\theta(u)}{\varphi(u)}\,\rmd u + \int_0^t\frac{\nu(u)}{\varphi(u)}\,\rmd B_u}
    \end{align*}
    is the unique strong solution (see Example 2.3 by \cite{Cetin_2016_Markov}) of \eqref{eq:tdOU_SDE}, where again $\varphi(t) = \exp\lrc{-\int_0^t \theta(u)\,\rmd u}$. That is, $X$ is a GM process with mean $m$ and covariance factorization $r_1$ and $r_2$ given by \eqref{eq:tdOU_mean_1} and \eqref{eq:tdOU_cov_fac}, respectively. 
    
    Relations \eqref{eq:cov_2} and \eqref{eq:cov_3} are trivial to check. From \eqref{eq:theta_explotion}, it follows \eqref{eq:cov_4}. The continuous differentiability of $\theta$, $\kappa$, and $\nu$ implies \eqref{eq:cov_5}. Using \eqref{eq:theta_sigma} and integrating by parts we get that 
    \begin{align}\label{eq:r1_theta}
        r_1(t) = 1 - \varphi(t).
    \end{align}
    Then, \eqref{eq:cov_6} holds, as $r_1(T) := \lim_{t\rightarrow T}r_1(t) = 1$ and $r_2(0) = 1$.
\end{proof}

\begin{remark}\label{rm:r2_and_beta_dec_r1_inc}
    After condition \eqref{eq:theta_sigma} and relation \eqref{eq:tdOUB_coeff_mean_cov}, we get that $r_2'(t)r_2(t) < 0$ for all $t\in(0, T)$. Hence, since $r_2$ is continuous and does not vanish in $[0, T)$, it can be chosen as either positive and decreasing, or negative and increasing. In \eqref{eq:dGMB_cov_fac}, the positive decreasing version is chosen, which is reflected by the fact that $\beta_T > 0$ is assumed in representation \eqref{eq:BM_rep}. Since $\beta_T = r_2$, then $\beta_T$ is also decreasing. Likewise, \eqref{eq:dGMB_cov_fac} and \eqref{eq:r1_theta} indicate that $r_1$ is chosen as positive and increasing.
\end{remark}

One could argue that defining a GMB should only require the process to degenerate at $t = 0$ and $t = T$, which is equivalent to \eqref{eq:cov_1}--\eqref{eq:cov_4}. GMBs defined in this way are not necessarily derived from conditioning a GM process, as it is assumed in representation \eqref{def:GMB_derived}. Indeed, consider the Gaussian process $X = \{X_t\}_{t\in[0, 1]}$ with zero mean and covariance function $R(t_1, t_2) = r_1(t_1\wedge t_2)r_2(t_1 \vee t_2)$ for all $t_1, t_2 \in [0, 1]$, where $r_1(t) = t^2(1-t)$ and $r_2(t) = t(1 - t)$. Lemma \ref{lm:GM_def} entails that $R$ is a valid covariance function and $X$ is Markovian. Moreover, since $r_1(0) = r_2(1) = 0$, $X$ is a bridge from $(0, 0)$ to $(1, 0)$. However, $r_1(0) = r_2(0) = 0$. That is, \eqref{eq:cov_6} fails and, hence, $X$ does not satisfy definition \eqref{def:GMB_cov}. Recognizing the differences between both definitions of GMBs, we adopt that in which a GM process is conditioned to take deterministic values at some initial and future time, since representation \eqref{eq:GMB_to_BM} is key to our results in Section~\ref{sec:solution}. It reveals the (linear) dependence of the mean with respect to $x$ and $z$, and it clarifies the relation between OUBs and GMBs in \eqref{def:GMB_OUB}.

Notice that a higher smoothness of the GMB mean and covariance factorization is assumed in all alternative characterizations in Proposition \ref{pr:GMB}. Clearly, this is a useful assumption to define GMBs, but not necessary. We discuss this in Remark \ref{rm:hypothesis}. In the rest of the paper, we implicitly assume the twice continuous differentiability of the mean and covariance factorization every time we mention a GMB.

Although easily obtainable from \eqref{eq:tdOUB_coeff_mean_cov}, for the sake of reference we write down the explicit relation between the BM representations \eqref{eq:GMB_to_BM} and the OUB representation \eqref{eq:tdOU_SDE}, namely:
\begin{align}\label{eq:tdOUB_coeff_BM}
\left\{
\begin{aligned}
    \theta(t) &= -\beta'_T(t)/\beta_T(t),
    \\
    \kappa(t) &= \alpha(t) - \beta_T(t)/\beta'_T(t) (\alpha'(t) + (z - \alpha(T))\beta_T(t)\gamma'_T(t)),
    \\
    \nu(t) &= \beta_T(t)\sqrt{\gamma'_T(t)}.
\end{aligned}
\right.
\end{align}

It is also worth mentioning that condition \eqref{eq:theta_sigma}, which is necessary and sufficient for an OU process to be an OUB, was also recently found in \cite[Theorem 3.1]{Hildebrandt_2020_pinned} for the case where $\kappa$ is assumed constant.


Finally, we rely on the classic OU process to illustrate the characterization in Lemma~\ref{lm:GM_def} and the connection between all alternative definitions in Proposition \ref{pr:GMB}.

\begin{example}[Ornstein--Uhlenbeck bridge]\ \\
    Let $\wt{X} = \lrc{X_t}_{t\in\R_+}$ be an OU process. That is, the unique strong solution of the SDE
    \begin{align*}
        \rmd X_t = aX_t\,\rmd t + c\,\rmd B_t,\quad t\in(0, T),
    \end{align*}
    where $\lrc{B_u}_{u\in\R_+}$ is a standard BM, and $a\in\R$, $c \in \R_+$. $\wt{X}$ is a time-continuous GM process non-degenerated on $[0, T]$. Its mean 
    and covariance factorization are twice continuously differentiable. In fact, they take the form
    \begin{align*}
        \wt{m}(t) &= \EEsp{\wt{X}_t} = \wt{X}_0e^{at}, \\
        \wt{R}(t_1, t_2) &= \Cov{\wt{X}_{t_1}}{\wt{X}_{t_2}} = \wt{r}_1(t_1\wedge t_2)\wt{r}_2(t_1\vee t_2), \\
        \wt{r}_1(t) &= \sinh(at), \quad \wt{r}_2(t) = c^2e^{at}/a.
    \end{align*}
    Note that $\wt{m}$, $\wt{r}_1$, and $\wt{r}_2$ satisfy conditions \ref{eq:cov_fac}--\ref{eq:cov_ratio} from Lemma \ref{lm:GM_def}.
    
    Let $\smash{X = \{X_u\}_{u\in [0, T]}}$ be a GM process defined on the same probability space as $\wt{X}$
    , for some $T > 0$. In agreement to Proposition \ref{pr:GMB}, the following statements are equivalent:
    \begin{enumerate}[label=(\textit{\roman{*}}), ref=\textit{\roman{*}}]
    \item $X$ results after conditioning $\wt{X}$ to $\wt{X}_0 = x$ and $\wt{X}_T = z$ in the sense of \eqref{pr:GMB_to_BM_1} from Proposition \ref{pr:GMB}, for some $x\in\R$ and $(T, z)\in\R_+\times\R$.

    \item The mean and covariance factorization of $X$ are twice continuously differentiable, and they satisfy conditions \eqref{eq:cov_1}--\eqref{eq:cov_6}. In fact, they take the form
    \begin{align*}
        m(t) &= \EEsp{X_t} =  (x\sinh(a(T - t)) + z\sinh(at))/\sinh(aT), \\ 
        R(t_1, t_2) &= \Cov{X_{t_1}}{X_{t_2}} = r_1(t_1\wedge t_2)r_2(t_1\vee t_2), \\
        r_1(t) &= \sinh(at)/\sinh(aT), \quad r_2(t) = c^2\sinh(a(T-t))/a,
    \end{align*}
    which follows after working out formulae \eqref{eq:dGMB_mean} and \eqref{eq:dGMB_cov_fac} (see also Proposition 3.3 from \cite{Barczy_2013_sample}).

    \item $X_T=z$ and, on $[0, T)$, $X$ admits the following representation
    \begin{align*}
        X_t 
        =\wt{X}_0e^{at} + \frac{c^2\sinh(a(T-t))}{a}\lrp{(z - \wt{X}_0e^{aT})\gamma_T(t) + \lrp{B_{\gamma_T(t)} + \frac{a(x - \wt{X}_0)}{c^2\sinh(aT)}}\!}\!.
    \end{align*}
    
    This expression does not depend on $\wt{X}_0$, indeed, after some manipulation it simplifies in 
    \begin{align*}
            X_t 
            = \frac{\sinh(a(T-t))}{\sinh(aT)}x + \frac{\sinh(at)}{\sinh(aT)}z + \frac{c^2\sinh(a(T-t))}{a}B_{\gamma_T(t)}\;,
    \end{align*}
    which is in alignment with the ``space-time transform'' representation in \cite{Barczy_2013_sample}.

    \item $X$ is the unique strong solution of the SDE
    \begin{align*}
        \rmd X_t = \theta(t)(\kappa(t) - X_t)\,\rmd t + \nu(t)\,\rmd B_t,\quad t\in(0, T),
    \end{align*}
    with initial condition $X_0 = x$ and
    \begin{align*}
    \left\{
    \begin{aligned}
        \theta(t) &= \coth(a(T-t)),
        \\
        \kappa(t) &= z/\cosh(a(T-t)),
        \\
        \nu(t) &= c.
    \end{aligned}
    \right.
    \end{align*}
    These expressions for the drift and volatility terms of $X$ come from \eqref{eq:tdOUB_coeff_BM}, and are in agreement with Equation (3.2) in \cite{Barczy_2013_sample}. Conditions \eqref{eq:theta_explotion} and \eqref{eq:theta_sigma} follow straightforwardly.
    \end{enumerate}
\end{example}

\section{Two equivalent formulations of the OSP}\label{sec:formulation}

For $0 \leq t < T$, let $X = \{X_u\}_{u\in[0, T]}$ be a real-valued, time-continuous GMB with $X_T = z$, for some $z\in\R$. Define the finite-horizon OSP
\begin{align}\label{eq:value_OSP_GMB}
    V_{T, z}(t, x) &:= \sup_{\tau\leq T - t }\EEs{X_{t + \tau}}{t, x},
\end{align}
where $V_{T, z}$ is the value function and $\EE_{t, x}$ is the mean operator with respect to the probability measure $\PPr_{t, x}$ such that $\PPr_{t, x}(X_t = x) = 1$. The supremum in \eqref{eq:value_OSP_GMB} is taken among all random times $\tau$ such that $t + \tau$ is a stopping time for $X$, although, for simplicity, we will refer to $\tau$ as a stopping time from now on.

Likewise, consider a BM $\{Y_u\}_{u\in\R_+}$ on the probability space $(\Omega, \cF, \Pr)$, and define the infinite-horizon OSP
\begin{align}\label{eq:value_OSP_BM}
    W_{T, z}(s, y) := \sup_{\sigma\geq 0}\Es{G_{T, z}\lrp{s + \sigma, Y_{s + \sigma}}}{s, y},
\end{align}
for $(s, y)\in\R_+\times\R$, where $\Pr_{s, y}$ and $\E_{s, y}$ have analogous definitions to those of $\PPr_{t, x}$ and $\EE_{t, x}$, that is, $Y_{s + u} = y + B_u$ under $\Pr_{s, y}$, where $\lrc{B_u}_{u\in\R_+}$ is a standard BM. The supremum in \eqref{eq:value_OSP_BM} is taken among the stopping times of $\{Y_{s+u}\}_{u\in\R_+}$, and the (gain) function $G_{T, z}$ takes the form
\begin{align}\label{eq:gain_BM}
    G_{T, z}(s, y) := \alpha(\gamma_T^{-1}(s)) + \beta_T(\gamma_T^{-1}(s))\lrp{(z - \alpha(T))s + y},
\end{align}
for $\alpha$, $\beta_T$, and $\gamma_T$ as in \eqref{pr:GMB_to_BM_1}--\eqref{pr:GMB_to_BM_3} from Proposition \ref{pr:GMB}.

Note that we have used different notations for the probability and expectation operators in the OSPs \eqref{eq:value_OSP_GMB} and \eqref{eq:value_OSP_BM}. The intention is to emphasize the difference between the probability spaces relative to the original GMB and the resulting BM. We shall keep this notation for the rest of the paper.

Solving \eqref{eq:value_OSP_GMB} and \eqref{eq:value_OSP_BM} means providing a tractable expression for $V_{T, z}(t, x)$ and $W_{T, z}(s, y)$, as well as finding stopping times (if they exist) $\tau^* = \tau^*(t, x)$ and $\sigma^* = \sigma^*(s, y)$ such that
\begin{align*}
    V_{T, z}(t, x) = \EEs{X_{t + \tau^*}}{t, x}, \quad
    W_{T, z}(s, y) = \Es{G_{T,z}\lrp{s + \sigma^*, Y_{s + \sigma^*}}}{s, y}.
\end{align*}
In such a case, $\tau^*$ and $\sigma^*$ are called Optimal Stopping Times (OSTs) for \eqref{eq:value_OSP_GMB} and \eqref{eq:value_OSP_BM}, respectively.

We claim that the OSPs \eqref{eq:value_OSP_GMB} and \eqref{eq:value_OSP_BM} are equivalent in the sense specified in the following proposition. In summary, the representation \eqref{eq:GMB_to_BM} equates the original GMB to the BM transformed by the gain function $G_{T, z}$, and \eqref{pr:GMB_to_BM_3} changes the finite horizon~$T$ into an infinite horizon.

\begin{proposition}[Equivalence of the OSPs]\label{pr:OSP_equivalence}\ \\
    Let $V$ and $W$ be the value functions in \eqref{eq:value_OSP_GMB} and \eqref{eq:value_OSP_BM}. For $(t, x)\in[0,T]\times\R$, let $s = \gamma_T(t)$ and $y = \lrp{x - \alpha(t)}/\beta_T(t) - \gamma_T(t)(z - \alpha(T))$. Then,
    \begin{align}\label{eq:equivalence_value}
        V_{T, z}(t, x) = W_{T, z}\lrp{s, y}.
    \end{align}
    Moreover, $\tau^* = \tau^*(t, x)$ is an OST for $V_{T, z}$ if and only if $\sigma^* = \sigma^*(s, y)$, defined such that $s + \sigma^* = \gamma_T(t + \tau^*)$, is an OST for $W$.
\end{proposition}

\begin{proof} 
    From \eqref{eq:GMB_to_BM}, we have the following representation for $X_{t+u}$ under $\PPr_{t, x}$:
    \begin{align*}
        X_{t+u} &= \alpha(t+u) + \beta(t+u)\lrp{(z-\alpha(T))\gamma_T(t+u) + \lrp{B_{\gamma_T(t+u)} + \frac{X_0 - \alpha(0)}{\beta(0)}}} \\
        &= G_{T, z}\lrp{\gamma_T(t + u), \lrp{B_{\gamma_T(t+u)} + \frac{X_0 - \alpha(0)}{\beta(0)}}} \\
        &= G_{T, z}\lrp{\gamma_T(t + u), \lrp{B_{\gamma_T(t+u)}-B_{\gamma_T(t)}+\frac{X_t - \alpha(t)}{\beta_T(t)} - (z - \alpha(t))\gamma_T(t)}},
    \end{align*}
    where, in the last equation, we used the relation
    \begin{align*}
        B_{\gamma_T(t)} + \frac{X_0 - \alpha(0)}{\beta_T(0)} = \frac{X_t - \alpha(t)}{\beta_T(t)} - (z - \alpha(t))\gamma_T(t).
    \end{align*}

    Let $Y_{s+v}:= B'_v+y$ and $B'_v:=B_{s+v}-B_s$, with $\{B'_v\}_{v\in\R_+}$ being a standard $\Pr_{s, y}$-BM. We recall that we use $\Pr$ instead of $\PPr$ to emphasize the time-space change, although the measure remains the same.
    
    We have that 
    \begin{align*}
        X_{t+u} 
        &= G_{T, z}\lrp{\gamma_T(t + u), Y_{\gamma_T(t+u)}} .
    \end{align*}
    
    For every stopping time $\tau$ of $\{X_{t + u}\}_{u\in[0, T-t]}$, consider the stopping time $\sigma$ of $\{Y_{s + u}\}_{u\in\R_+}$ such that $s + \sigma = \gamma_T(t + \tau)$. 
    Hence, \eqref{eq:equivalence_value} follows from the following sequence of equalities:
    \begin{align*}
        V_{T, z}(t, x) &= \sup_{\tau \leq T-t} \EEs{X_{t + \tau}}{t, x} 
        = \sup_{\sigma \geq 0} \Es{G_{T, z}\lrp{s + \sigma, Y_{s + \sigma}}}{s, y} 
        = W_{T, z}\lrp{s, y}.
    \end{align*}
    Furthermore, suppose that $\tau^* = \tau^*(t, x)$ is an OST for \eqref{eq:value_OSP_GMB} and that there exists a stopping time $\sigma' = \sigma'(s, y)$ that performs better than $\sigma^* = \sigma^*(s, y)$ in \eqref{eq:value_OSP_BM}. Consider $\tau' = \tau'(t, x)$ such that $t + \tau' = \gamma_T^{-1}(s + \sigma')$. Then,
    \begin{align*}
        \EEs{X_{t + \tau'}}{t, x} = \Es{G_{t, T}(s + \sigma', Y_{s + \sigma'})}{s, y} > \Es{G_{t, T}(s + \sigma^*, Y_{s + \sigma*})}{s, y} = \EEs{X_{t + \tau^*}}{t, x},
    \end{align*}
    which contradicts the fact that $\tau^*$ is optimal. Using similar arguments, we can obtain the reverse implication, that is, if $\sigma^*$ is an OST for \eqref{eq:value_OSP_BM}, then $\tau^*$ is an OST for \eqref{eq:value_OSP_GMB}.
\end{proof}

\section{Solution of the infinite-horizon OSP}\label{sec:solution}

We have shown that solving \eqref{eq:value_OSP_GMB} is equivalent to solving \eqref{eq:value_OSP_BM}, which is expressed in terms of a simpler BM. In this section we leverage that advantage to solve \eqref{eq:value_OSP_BM} but, first, we rewrite it with a cleaner notation that hides its explicit connection with the original OSP, and allows us to treat \eqref{eq:value_OSP_BM} as a standalone problem. 

Let $\lrc{Y_u}_{u\in\R_+}$ be a BM on the probability space $\lrp{\Omega, \cF, \Pr}$. Define the probability measure $\Pr_{s,y}$ such that $\Pr_{s, y}(Y_s = y) = 1$.
Consider the OSP
\begin{align}\label{eq:value_OSP}
    W(s, y) := \sup_{\sigma \geq 0}\Es{G(s + \sigma, Y_{s+\sigma})}{s, y} = \sup_{\sigma \geq 0}\Esp{G(s + \sigma, Y_{\sigma} + y)},
\end{align}
where $\E$ and $\E_{s, y}$ are the mean operators with respect to $\Pr$ and $\Pr_{s, y}$, respectively. The supremum in \eqref{eq:value_OSP} is taken among the stopping times of $Y = \lrc{Y_{s+u}}_{u\in\R_+}$. The (gain) function $G$ takes the form
\begin{align}\label{eq:gain_function}
    G(s, y) = a_1(s) + a_2(s)\lrp{c_0s + y},
\end{align}
where $c_0\in\R$ and $a_1, a_2:\R_+\rightarrow\R$ are assumed to be such that:
\begin{subequations}
\begin{align}
    a_1 \text{ and } a_2 \text{ are twice continuously differentiable}, \label{assumption:continuous_diff} \\
    \text{$a_1$, $a_1'$, $a_1''$, $a_2$, $a_2'$, and $a_2''$ are bounded}; \label{assumption:bounded} \\
    \text{there exists } c_1\in\R \text{ such that} \lim_{s\rightarrow\infty} a_1(s) = c_1; \label{assumption:limit_1} \\
    \text{for all } s\in\R,\ a_2(s) > 0; \label{assumption:a_2>0} \\
    \text{there exists } c_2\in\R \text{ such that} \lim_{s\rightarrow\infty} a_2(s)s = c_2; \label{assumption:limit_2} \\
    \text{for all } s\in\R,\ a_2'(s) < 0. \label{assumption:a_2'<0}
\end{align}
\end{subequations}

Assumptions \eqref{assumption:continuous_diff}--\eqref{assumption:a_2'<0} do not further restrict the class of GMBs considered in Proposition \ref{pr:GMB}. Indeed, \eqref{assumption:continuous_diff}--\eqref{assumption:bounded} are implied by the twice continuous differentiability of the GMB's mean and covariance factorization, while \eqref{assumption:limit_1}--\eqref{assumption:a_2'<0} are obtained from the degenerative nature of the GMB. In fact, the infinite-horizon OSP \eqref{eq:value_OSP} under assumptions \eqref{assumption:continuous_diff}--\eqref{assumption:a_2'<0} is equivalent to the finite-horizon OSP \eqref{eq:value_OSP_GMB} with a GMB as the underlying process. The following remarks shed light on this equivalence.

\begin{remark}\label{rm:c0-a1-a2}
    Equation \eqref{eq:gain_function}, as well as assumptions \eqref{assumption:limit_1}--\eqref{assumption:limit_2}, come after \eqref{eq:gain_BM} and \eqref{pr:GMB_to_BM_1}--\eqref{pr:GMB_to_BM_3} from Proposition \ref{pr:GMB}. Indeed, the constant $c_0$ and the functions $a_1$ and $a_2$ are taken such that $c_0 = z - \alpha(T)$, $a_1(s) = \alpha(\gamma_{T}^{-1}(s))$, and $a_2(s) = \beta_T(\gamma_T^{-1}(s))$. 
\end{remark}    

\begin{remark}\label{rm:hypothesis}
    Assumptions \eqref{assumption:continuous_diff} and \eqref{assumption:bounded} are derived from the twice continuous differentiability of $\alpha$, $\beta_T$, and $\gamma_T$. 
    These assumptions are used to prove smoothness properties of the value function and the OSB. The assumptions on the first derivatives are used to prove the Lipschitz continuity of the value function (see Proposition \ref{pr:value_Lipschitz}), while the ones on the second derivatives are required to prove the local Lipschitz continuity of the OSB (see Proposition \ref{pr:OSB_Lipschitz}). 
\end{remark}

\begin{remark}
    The following relation, which we use recurrently throughout the paper, comes after \eqref{assumption:continuous_diff}, \eqref{assumption:bounded}, and \eqref{assumption:limit_2}:
    \begin{align}\label{limit_3}
        \lim_{s\rightarrow\infty} a_2'(s)s = 0.
    \end{align}
    Alternatively, \eqref{limit_3} can be directly derived from \eqref{eq:dGMB_cov_fac} and the fact that $\lim_{s\rightarrow\infty}a_2(s) = 0$. Indeed,
    \begin{align*}
        \lim_{s\rightarrow\infty}a_2'(s)s =&\; \lim_{s\rightarrow\infty}a_2'(s)s + a_2(s) = \lim_{s\rightarrow\infty}\partial_s \lrb{a_2(s)s} = \lim_{s\rightarrow\infty} \partial_s r_1(\gamma_T^{-1}(s)) = \lim_{t\rightarrow T} \frac{r_1'(t)}{\gamma_T'(t)} \\
        =&\; \lim_{t\rightarrow T}\frac{r_1'(t)r_2^2(t)}{r_1'(t)r_2(t) - r_1(t)r_2'(t)} = 0,
    \end{align*}
    where $\partial_s$ denotes the derivative with respect to the variable $s\in\R_+$. In the last equality we used that $0 \leq r_1'(t)/r_2'(t) \leq r_1(t)/r_2(t)$, which comes after $r_1$ and $r_2$ being, respectively, an increasing and a decreasing function (see Remark \ref{rm:r2_and_beta_dec_r1_inc}).
    
    Likewise, \eqref{limit_3} along with the L'Hôpital rule implies that
    \begin{align}\label{limit_4}
        \lim_{s\rightarrow\infty} a_2''(s)s^2 = -\lim_{s\rightarrow\infty} a_2'(s)s = 0.
    \end{align}
    Again, \eqref{limit_4} can be obtained from its representation in terms of the covariance factorization given by $r_1$ and $r_2$,
    \begin{align*}
        \lim_{s\rightarrow\infty}a_2''(s)s^2 =&\;
        \lim_{s\rightarrow\infty}\partial_{ss} \lrb{a_2(s)s}s = \lim_{s\rightarrow\infty} \partial_{ss} r_1(\gamma_T^{-1}(s))\gamma_T(\gamma_T^{-1}(s)) \\
        =&\; \lim_{s\rightarrow\infty} \partial_{s} \frac{r_1'(\gamma_T^{-1}(s))}{\gamma_T'(\gamma_T^{-1}(s))}\gamma_T(\gamma_T^{-1}(s)) 
        = \lim_{t\rightarrow T} \lrp{\frac{r_1''(t)}{(\gamma_T'(t))^2} - \frac{r_1(t)\gamma_T''(t)}{(\gamma_T'(t))^3}}\gamma_T(t) \\
        =&\;\lim_{t\rightarrow T} \lrp{\frac{r_1^2(t)r_2^3(t)\blrp{r_1'(t)r_2''(t) - r_1''(t)r_2'(t)}}{\blrp{\!r_1'(t)r_2(t) - r_1(t)r_2^3(t)}} -  \frac{2r_1(t)r_2^2(t)\blrp{r_1'(t)}^2}{\blrp{\!r_1'(t)r_2(t) - r_1(t)r_2'(t)}^2}} \\
        =&\; 0,
    \end{align*}
    where $\partial_{ss}$ indicates the second derivative with respect to $s$.
\end{remark}
    
\begin{remark}
    Assumption \eqref{assumption:a_2'<0} is needed to derive the boundedness of the OSB (see Proposition \ref{pr:OSB_shape} and Remark \ref{rmk:assumption_a_2'<0}). Similarly to Assumptions \eqref{assumption:continuous_diff}--\eqref{assumption:limit_2}, Assumption \eqref{assumption:a_2'<0} can be obtained from the regularity of the underlying GMB already used in Section \ref{sec:GMB}, and does not impose any further restrictions. Specifically, Assumption \eqref{assumption:a_2'<0} is equivalent to 
    condition $\theta(t) > 0$ for all $t\in[0, T]$, in the OUB representation \eqref{def:GMB_OUB} from Proposition \ref{pr:GMB}, and to $\beta_T(t) = r_2(t) > 0$ and $\beta_T'(t) = r_2'(t)< 0$, in terms of representations \eqref{def:GMB_BM} and \eqref{def:GMB_cov} (see Remark \ref{rm:r2_and_beta_dec_r1_inc}).
\end{remark}

Notice that \eqref{assumption:limit_1}, \eqref{assumption:limit_2}, and \eqref{limit_3}, together with the law of the iterated logarithm, imply that, for all $(s, y)\in\R_+\times\R$,
\begin{align}\label{eq:gain_convergence}
    \Pr_{s,y}\text{-\!\!}\lim_{u\to\infty}G(s + u, Y_{s+u}) = c_1 + c_0c_2.
\end{align}

For later reference, let us introduce the notation
\begin{align}\label{eq:bounds}
\left.
    \begin{aligned}
        A_1 &:= \sup_{s\in\R_+} \lrav{a_1(s)},\quad A_1' := \sup_{s\in\R_+} \lrav{a_1'(s)},\quad A_1'' := \sup_{s\in\R_+} \lrav{a_1''(s)},  \\
        A_2 &:= \sup_{s\in\R_+} \lrav{a_2(s)},\quad A_2' := \sup_{s\in\R_+} \lrav{a_2'(s)},\quad A_2'' := \sup_{s\in\R_+} \lrav{a_2''(s)}, \\
        A_3 &:= \sup_{s\in\R_+} \lrav{a_2(s)s},\hspace{0.24cm} A_3' := \sup_{s\in\R_+} \lrav{a_2'(s)s}, \hspace{0.18cm} A_3'' := \sup_{s\in\R_+} \lrav{a_2''(s)s}.
    \end{aligned}
\right\}
\end{align}
In addition, we will often require the expression of the partial derivatives of $G$, namely,
\begin{align} 
    \partial_t G(s, y) &= a_1'(s) + c_0a_2(s) + a_2'(s)(c_0s + y), \label{eq:G_t} \\
    \partial_x G(s, y) &= a_2(s). \label{eq:G_y}
\end{align}
Here and thereafter, $\partial_t$ and $\partial_x$ stand, respectively, for the differential operator with respect to time and space.

Notice that \eqref{assumption:limit_2} guarantees the existence of $m > 0$ such that $\lrav{a_2(s)} \leq (1 + m)/s$ for all $s\geq 1$, which, combined with the boundedness of $a_1$, $a_2$, and $s\mapsto a_2(s)s$, entails the following bound with $A = \max\{A_1 + |c_0|A_3, A_2\}$:
\begin{align}
    \E_{s, y}\bigg[\!&\sup_{u\in\R_+} \lrav{G\lrp{s + u, Y_{s+u}}}\bigg]\nonumber\\
    \leq&\; \sup_{u\in\R_+}\lrav{a_1(u) + a_2(u)(c_0u + y)} + \Espbigg{\sup_{u\in\R_+} \lrav{a_2(s + u)Y_u}} \nonumber \\
    \leq&\; A(1 + |y|) + \Espbigg{\sup_{u\in\R_+}\lrav{a_2(s + u)Y_u}} \nonumber \\
    \leq&\; A(1 + |y|) +\!\max_{u\leq 1\vee(1 - s)}\!\lrav{a_2(s + u)}\Esp{\sup_{u\leq 1\vee(1 - s)}\!\lrav{Y_u}}\! + \Espbigg{\sup_{u\geq 1\vee(1 - s)}\!\lrav{a_2(s + u)Y_u}} \nonumber \\
    \leq&\; A(1 + |y|) + \max_{u\leq 1}\lrav{a_2(u)}\Esp{\sup_{u\leq 1}\lrav{Y_u}} + (1 + m)\Espbigg{\sup_{u\geq 1}\lrav{Y_u} \big/ u} \nonumber \\
    =&\; A\lrp{1 + \lrp{\lrav{y} + \Esp{\sup_{u\leq 1}\lrav{Y_u}}}} + (1 + m)\Esp{\sup_{u\geq 1}\lrav{Y_{1/u}}} \nonumber \\
    =&\; A\lrp{1 + \lrp{\lrav{y} + \Esp{\sup_{u\leq 1}\lrav{Y_u}}}} + (1 + m)\Esp{\sup_{u\leq 1}\lrav{Y_u}} < \infty. \label{eq:dominated_convergence_condition}
\end{align}
In the last equality, the time-inversion property of the BM was used.

The continuity of $G$ alongside \eqref{eq:dominated_convergence_condition} implies the continuity of $W$. However, given the assumptions \eqref{assumption:continuous_diff}--\eqref{assumption:limit_2}, one can obtain higher smoothness for the value function, namely its Lipschitz continuity, as shown in the proposition below.

\begin{proposition}[Lipschitz continuity of the value function]\label{pr:value_Lipschitz}\ \\
    For any bounded set $\cR\subset \mathbb{R}$ there exists $L_\cR > 0$ such that
    \begin{align}\label{eq:value_Lipschitz}
        \lrav{W(s_1, y_2) - W(s_2, y_2)} \leq L_\cR(|s_1 - s_2| + |y_1 - y_2|),
    \end{align}
    for all $(s_1, y_1), (s_2, y_2)\in\R_+\times\cR$.
\end{proposition}

\begin{proof}
    For any $(s_1, y_1), (s_2, y_2)\in\R_+\times \cR$, the following equality holds:
    \begin{align*}
        W(s_1, y_1) - W(s_2, y_2) =&\; \sup_{\sigma\geq 0}\Es{G\lrp{s_1 + \sigma, Y_{s_1 + \sigma}}}{s_1, y_1} - \sup_{\sigma\geq 0}\Es{G\lrp{s_1 + \sigma, Y_{s_1 + \sigma}}}{s_1, y_2} \\
        &+ \sup_{\sigma\geq 0}\Es{G\lrp{s_1 + \sigma, Y_{s_1 + \sigma}}}{s_1, y_2} - \sup_{\sigma\geq 0}\Es{G\lrp{s_2 + \sigma, Y_{s_2 + \sigma}}}{s_2, y_2}.
    \end{align*}
    Since $|\sup_\sigma a_\sigma - \sup_\sigma b_\sigma|\leq \sup_\sigma|a_\sigma - b_\sigma|$, and due to Jensen's inequality,
    \begin{align}
        \bigg|\sup_{\sigma\geq 0}\E_{s_1,y_1}&\, [G\lrp{s_1 + \sigma, Y_{s_1 + \sigma}}] - \sup_{\sigma\geq 0}\Es{G\lrp{s_1 + \sigma, Y_{s_1 + \sigma}}}{s_1, y_2}\bigg| \nonumber \\
        &\leq \Esp{\sup_{u \geq 0}\lrav{G\lrp{s_1 + u, Y_u + y_1} - G\lrp{s_1 + u, Y_u + y_2}}} \nonumber \\
        &= \sup_{u \geq 0}\lrav{a_2(s_1 + u)(y_1 - y_2)} \nonumber \\
        &\leq A_2\lrav{y_1 - y_2}.
        \label{eq:value_Lipschitz_space}
    \end{align}
    Likewise,
    \begin{align}
         \bigg|\sup_{\sigma\geq 0}\E_{s_1,y_2}&\,G\lrp{s_1 + \sigma, Y_{s_1 + \sigma}}] - \sup_{\sigma\geq 0}\Es{G\lrp{s_2 + \sigma, Y_{s_2 + \sigma}}}{s_2, y_2}\bigg| \nonumber \\
        &\leq \Esp{\sup_{u \geq 0}\lrav{G\lrp{s_1 + u, Y_u + y_2} - G\lrp{s_2 + u, Y_u + y_2}}} \nonumber \\
        &= \Esp{\sup_{u \geq 0}\lrav{\partial_t G\lrp{\eta_u, Y_u + y_2}(s_1 - s_2)}} \nonumber \\
        &\leq \lrp{A_1' + (A_3' + A_2)|c_0| + A_2'\lrp{\sup_{y\in\cR}\lrc{y} + \Esp{\sup_{u \geq 0}\lrav{Y_u}}}}|s_1 - s_2|, \label{eq:value_Lipschitz_time}
    \end{align}
    where $\eta_u\in(s_1\wedge s_2 + u, s_1\vee s_2 + u)$ comes from the mean value theorem, which, along with \eqref{eq:G_t}, was used to derive the last inequality. Constants $A_1'$, $A_2$, $A_2'$, and $A_3'$ were defined in \eqref{eq:bounds}. We finally get \eqref{eq:value_Lipschitz} after merging \eqref{eq:value_Lipschitz_space} and \eqref{eq:value_Lipschitz_time}.
\end{proof}

Define $\sigma^* = \sigma^*(s, y) := \inf\lrc{u\in\R_+ : (s + u, Y_{s+u}) \in \cD}$, where the closed set
\begin{align*}
    \cD := \lrc{(s, y)\in\R_+\times\R: W(s, y) = G(s, y)},
\end{align*}
is called the \emph{stopping set}. The continuity of $W$ and $G$ (it suffices lower semi-continuity of $W$ and upper semi-continuity of $G$) along with \eqref{eq:dominated_convergence_condition} and \eqref{eq:gain_convergence}, guarantees that $\sigma^*$ is an OST for \eqref{eq:value_OSP} (see Corollary 2.9 and Remark 2.10 in \cite{Peskir_2006_optimal}), meaning that
\begin{align}\label{eq:value_OST}
    W(s, y) = \Es{G(s + \sigma^*, Y_{s + \sigma^*})}{s, y}.
\end{align}
Applying Itô's lemma to \eqref{eq:value_OSP} and \eqref{eq:value_OST}, we get a martingale term $\int_0^u a_2(s + r)\,\mathrm{d} B_r$ that results uniformly integrable as 
$\int_0^\infty a_2^2(s + r)\,\mathrm{d} r < \infty$, due to \eqref{assumption:limit_2}. Taking $\Pr_{s,y}$-expectation, this term vanishes and we get the following alternative representations of $W$:
\begin{align}
    W(s, y) - G(s, y) &= \sup_{\sigma\geq 0}\Es{\int_0^\sigma \InfGen G\lrp{s + u, Y_{s+u}}\,\rmd u}{s, y}\nonumber\\
    &= \Esbigg{\int_0^{\sigma^*(s, y)} \InfGen G\lrp{s + u, Y_{s+u}}\,\rmd u}{s, y},\label{eq:value_Ito}
\end{align}
where $\InfGen := \partial_t + \frac{1}{2}\partial_{xx}$ is the infinitesimal generator of the process $\lrc{\lrp{s, Y_s}}_{s\in \R_+}$ and the operator $\partial_{xx}$ is a shorthand for $\partial_x\partial_x$. Note that $\InfGen G = \partial_t G$.

Denote by $\cC$ the complement of $\cD$,
\begin{align*}
    \cC := \lrc{(s, y) \in\R_+\times\R : W(s, y) > G(s, y)},
\end{align*}
which is called the \emph{continuation set}. The boundary between $\cD$ and $\cC$ is the OSB and it determines the OST $\sigma^*$. 

In addition to the Lipschitz continuity, higher smoothness of the value function is achieved away from the OSB, as stated in the next proposition. We also determine the connection between the OSP \eqref{eq:value_OSP} and its associated free-boundary problem. For further details on this connection in a more general setting we refer to Section 7 of \cite{Peskir_2006_optimal}. 

\begin{proposition}[Higher smoothness of the value function and the free-boundary problem]\label{pr:value_smoothness}\ \\
    $W\in C^{1, 2}(\cC)$, that is, the functions $\partial_t W$, $\partial_x W$, and $\partial_{xx} W$ exist and are continuous on $\cC$. Additionally, $y\mapsto W(s, y)$ is convex for all $s\in\R_+$ and $\InfGen W = 0$ on $\cC$.
\end{proposition}

\begin{proof}
	The convexity of $W$ with respect to the space coordinate is a straightforward consequence of the linearity of $Y_{s+u}$ with respect to $y$ under $\Pr_{s, y}$, for all $s\in\R_+$. Indeed, it follows from \eqref{eq:value_OSP} that $W(s, ry_1 + (1-r)y_2) \leq rW(s, y_1) + (1-r)W(s, y_2)$, for all $y_1,y_2\in\R$ and $r\in[0, 1]$.
	
	Since $W$ is continuous on $\cC$ (see Proposition \ref{pr:value_Lipschitz}) and the coefficients in the parabolic operator $\InfGen$ are smooth enough (it suffices to require local $\alpha$-H\"older continuity), then standard theory of parabolic partial differential equations \cite[Section 3, Theorem 9]{Friedman_1964_partial} guarantees that, for an open rectangle $\cR\subset \cC$, the initial-boundary value problem
	\begin{align}\label{eq:PDE}
	    \left\{
	\begin{aligned}
		\InfGen f &= 0 && \text{in } \cR,  \\
    	f &= W && \text{on } \partial \cR,
	\end{aligned}
	    \right.
	\end{align}
	where $\partial \cR$ refers to the boundary of $\cR$, has a unique solution $f\in C^{1, 2}(\cR)$. Therefore, we can use Itô's formula on $f(s + u, Y_{s+u})$ at $u = \sigma_{\cR}$, that is, the first time $(s+u, Y_{s+u})$ exits $\cR$, and then take $\Pr_{s, y}$-expectation with $(s, y) \in \cR$, which guarantees the vanishing of the martingale term and yields, along with \eqref{eq:PDE} and the strong Markov property, the equalities $W(s, y) = \Es{W(s + \sigma_{\cR}, Y_{s + \sigma_{\cR}})}{s, y} = f(s, y)$.
	Since $W = G$ on $\cD$, it follows that $W\in C^{1, 2}(\cD)$.
\end{proof}

In addition to the partial differentiability of $W$, it is possible to provide relatively explicit forms for $\partial_t W$ and $\partial_x W$ by relying on representation \eqref{eq:value_Ito} and the fact that $a_1$ and $a_2$ are differentiable functions.

\begin{proposition}[Partial derivatives of the value function]\label{pr:W_t_&_W_x}\ \\
    For any $(s, y)\in\cC$, consider the OST $\sigma^* = \sigma^*(s, y)$. Then, 
    \begin{align}\label{eq:W_t}
        \hspace*{-0.2cm}\partial_t W(s, y) = \partial_t G(s, y) + \Esbigg{\int_s^{s + \sigma^*}\lrp{a_1''(u) + 2c_0a_2'(u) + a_2''(u)(c_0u + Y_u)}\, \rmd u}{s, y}
    \end{align}
    and
    \begin{align}\label{eq:W_x}
        \partial_x W(s, y) = \Es{a_2(s + \sigma^*)}{s, y}.
    \end{align}
\end{proposition}

\begin{proof}
    Since $\sigma^* = \sigma^*(s, y)$ is suboptimal for any initial condition other than $(s, y)$, then
    \begin{align*}
        \varepsilon^{-1}\lrp{W(s, y) - W(s - \varepsilon, y)} 
        &\leq \varepsilon^{-1}\Esp{G(s + \sigma^*, Y_{\sigma^*} + y) - G(s - \varepsilon + \sigma^*, Y_{\sigma^*} + y)},
    \end{align*}
    for any $0 < \varepsilon \leq s$. Hence, by letting $\varepsilon\rightarrow 0$ and recalling that $W\in C^{1,2}(\cC)$ (see Proposition \ref{pr:value_smoothness}), we get that, for $(s, y)\in\cC$,
    \begin{align}\label{eq:W_t<}
        \hspace*{-0.2cm}\partial_t W(s, y) &\leq \Es{\partial_t G(s + \sigma^*, Y_{s + \sigma^*})}{s, y} = \partial_t G(s, y) + \Esbigg{\int_0^{\sigma^*}\InfGen\partial_t G(s + u, Y_{s+u})\, \rmd u}{s, y}.
    \end{align}
    In the same fashion, we obtain that 
    \begin{align*}
        \varepsilon^{-1}\lrp{W(s + \varepsilon, y) - W(s, y)} 
        &\geq \varepsilon^{-1}\Esp{G(s + \varepsilon + \sigma^*, Y_{\sigma^*} + y) - G(s + \sigma^*, Y_{\sigma^*} + y)},
    \end{align*}
    which, after letting $\varepsilon\rightarrow 0$, yields \eqref{eq:W_t<} in the reverse direction. 
    Therefore, \eqref{eq:W_t} is proved after computing $\InfGen\partial_t G(s + u, Y_{s+u}) = \partial_{tt} G(s + u, Y_{s+u})$.
    
    To get the analog result for the space coordinate, notice that
    \begin{align*}
        \varepsilon^{-1}\lrp{W(s, y) - W(s, y - \varepsilon)} &\leq \varepsilon^{-1}\Esp{W(s + \sigma^*, Y_{\sigma^*} + y) - W(s + \sigma^*, Y_{\sigma^*} + y - \varepsilon)} \\
        &\leq \varepsilon^{-1}\Esp{G(s + \sigma^*, Y_{\sigma^*} + y) - G(s + \sigma^*, Y_{\sigma^*} + y - \varepsilon)} \\
        &= \Es{a_2(s + \sigma^*)}{s, y},
    \end{align*}
    while the same reasoning yields the inequality 
    \begin{align*}
        \varepsilon^{-1}\lrp{W(s, y + \varepsilon) - W(s, y)} \geq \Es{a_2(s + \sigma^*)}{s, y},
    \end{align*}
    and then, by letting $\varepsilon\rightarrow 0$, \eqref{eq:W_x} follows.
\end{proof}

Besides the regularity of the value function, that of the OSB is also key to solving the OSP. However, defined as the boundary between $\cD$ and $\cC$, the OSB admits little space for technical manipulations. The next proposition gives a handle on the OSB by showing that it is the graph of a bounded function of time, above which $\cD$ lies.

\begin{proposition}[Shape of the OSB]\label{pr:OSB_shape}\ \\
    There exists a function $b:\R_+\rightarrow\R$ such that 
    \begin{align*}
        \cD = \lrc{(s, y)\in\R_+\times\R : y \geq b(s)}.
    \end{align*}
    Moreover, $g(s) < b(s) < \infty$ for all $s\in\R_+$, where $g(s) := (-a_1'(s) - c_0(a_2(s) + a_2'(s)s))/a_2'(s)$.
\end{proposition}

\begin{proof}
    Define $b$ as
    \begin{align}\label{eq:OSB}
        b(s) := \inf\lrc{y: (s, y) \in \cD},\quad s\in\R_+.
    \end{align}
    The claimed shape for the stopping set 
    is a straightforward consequence of the decreasing behavior of $y\mapsto (W-G)(s, y)$ for all $s\in\R_+$, which comes after \eqref{assumption:a_2'<0}, \eqref{eq:G_t}, and~\eqref{eq:value_Ito}.
    
    To derive the lower bound for $b$, notice that, for all $(s, y)$ such that $\partial_t G(s, y) > 0$, we can pick a ball $\cB$ such that $(s, y)\in\cB$ and $\partial_t G > 0$ on $\cB$. After recalling \eqref{eq:value_Ito} and by letting $\sigma_\cB = \sigma_\cB(s, y)$ to be the first exit time of $Y^{s, y}$ from $\cB$, we get that
    \begin{align*}
    W(s, y) - G(s, y) \geq \Es{\int_0^{\sigma_\cB} \partial_t G\lrp{s + u, Y_{s+u}}\,\rmd u}{s, y} > 0,
    \end{align*}
    which means that $(s, y)\in\cC$. Finally, the claimed lower bound for $b$ comes after using \eqref{eq:G_t} and \eqref{assumption:a_2'<0} to realize that $\partial_t G(s, y) > 0$ if and only if $y < g(s)$. 
    
    We now prove that $b(s) < \infty$ for all $s \in\R_+$.
    Let $X = \big\{X_t\big\}_{t \in [0, T]}$ be the OUB representation of the process $s\mapsto G(s, Y_s)$, that is, 
    the unique strong solution of \eqref{eq:tdOU_SDE}, 
    with drift $\mu(t, x) = \theta(t)(\kappa(t) - x)$ and volatility (function) $\nu$. This GMB $X$ is well defined, as we can trace back functions $\alpha$, $\beta_T$, and $\gamma_T$ and values $T$ and $z$, such that the OSP \eqref{eq:value_OSP_BM} is in the form $\eqref{eq:value_OSP}$ (see Remark \ref{rm:c0-a1-a2}).
    
    In addition to $X$, define the OUBs $X^{(i)}$, for $i = 1, 2$, with volatility $\nu$ and drifts
    \begin{align*}
        \mu^{(1)}(t, x) = \theta(t)(K - x),\quad 
        \mu^{(2)}(t, x) = \frac{\ul{\nu}}{\ol{\nu}(T-t)}(K - x),
    \end{align*}
    respectively, where $K := \max\{\kappa(t) : t\in[0, T]\}$, $\ol{\nu} := \max\{\nu(t) : t\in[0, T]\}$, and $\ul{\nu} := \min\{\nu(t) : t\in[0, T]\}$. Consider the OSPs
    \begin{align*}
        V^{(0)}(t, x) &:= \sup_{\tau\leq T-t}\Es{X_{t+\tau}}{t, x}, \\
        V^{(1)}(t, x) &:= \sup_{\tau\leq T-t}\Es{X_{t+\tau}^{(1)}}{t, x}, \\
        V_K^{(2)}(t, x) &:= \sup_{\tau\leq T-t}\Es{K + |X_{t+\tau}^{(2)} - K|}{t, x},
    \end{align*}
    alongside their respective stopping sets $\cD^{(0)}$, $\cD^{(1)}$, and $\cD^{(2)}_K$.
    
    Notice that $\mu(t, x) \leq \mu^{(1)}(t, x)$ for all $(t, x)\in[0, T)\times\R$. Hence, $X_{t+u}\leq X_{t+u}^{(1)}$\ $\Pr_{t, x}$\nobreakdash-a.s. for all $u\in[0, T-t]$, as Corollary 3.1 in \cite{Peng_2006_necessary} states. This implies that $\cD^{(1)} \subset \cD^{(0)}$.
    
    On the other hand, it follows from \eqref{eq:theta_sigma} that $\theta(t) \geq \ul{\nu}/(\ol{\nu}(T-t))$, meaning that $\mu(t, x) \leq \mu^{(2)}(t, x)$ if and only if $x \geq K$. By using the same comparison result in \cite{Peng_2006_necessary}, we get the second inequality in the following sequence of relations:
    \begin{align*}
        X_{t+u}^{(1)} \leq K + |X_{t+u}^{(1)} - K| \leq K + |X_{t+u}^{(2)} - K|
    \end{align*}
    $\Pr_{t, x}$-a.s. for all $u\in[0, T-t]$. Hence, for a pair $(t, x)\in \cD_K^{(2)}$, we get that $V^{(0)}(t, x) \leq V^{(2)}_K(t, x) = x$, that is, $(t, x) \in \cD^{(1)}$ and, therefore, $\cD^{(2)}_K \subset \cD^{(0)}$. The OSP related to $V^{(2)}_K$ can be shown to account for a finite OSB. Specifically, it is a multiple of that of a BB (see \cite[Section 5]{DAuria_2020_class}). Then, $\cD^{(0)} \cap \lrp{\{t\}\times\R}$ is non-empty for all $t\in[0,  T)$, and the equivalence result in Proposition \ref{pr:OSP_equivalence} guarantees that so are the sets of the form $\cD \cap \lrp{\{t\}\times\R}$, meaning that the OSB $b$ is bounded from above.
\end{proof}

\begin{remark}\label{rmk:assumption_a_2'<0}
    Note that the same reasoning we used to derive the lower bound of $b$ in the proof of Proposition \ref{pr:OSB_shape} also implies that, if $a_2'(s) > 0$ for some $s\in\R_+$, then $(s, y)\in\cC$ for all $y > (-a_1'(s) - c(a_2(s) + a_2'(s)s))/a_2'(s)$, meaning that $b(s) = \infty$. To avoid this explosion of the OSB we impose $a_2'(s) < 0$ for all $s\in\R_+$ in \eqref{assumption:a_2'<0}.
\end{remark}

Summarizing, we have proved that $W$ satisfies the free-boundary problem
\begin{align*}
    \InfGen W(s, y) &= 0 &&\hspace*{-3.25cm} \text{for } y < b(t), 
    \\
    W(s, y) &> G(s, y) &&\hspace*{-3.25cm} \text{for } y < b(t),
    \\
    W(s, y) &= G(s, y) &&\hspace*{-3.25cm} \text{for } y\geq b(t).
\end{align*}
Since $b$ is unknown, an additional condition, generally known as the \emph{smooth-fit condition}, is needed to guarantee the uniqueness of the solution of this free-boundary problem. When $b$ is regular enough, it comes by making the value and the gain function coincide smoothly at the free boundary. 

The works of \cite{DeAngelis_2015_note}, \cite{Peskir_2019_continuity}, and \cite{DeAngelis_2019_Lipschitz} address the smoothness of the free boundary. For one-dimensional, time-homogeneous processes with locally Lipschitz-continuous drift and volatility, \cite{DeAngelis_2015_note} provides the continuity of the free boundary. \cite{Peskir_2019_continuity} works with the two-dimensional case in a fairly general setting, proving the impossibility of first-type discontinuities (second-type discontinuities are not addressed). \cite{DeAngelis_2019_Lipschitz} goes further by proving the local Lipschitz continuity of the free boundary in a higher-dimensional framework. In particular, local Lipschitz continuity suffices for the smooth-fit condition to hold (see Proposition \ref{pr:smooth-fit} ahead), 
which is the main drive to tailor the method of \cite{DeAngelis_2019_Lipschitz} to fit our settings in the next proposition. Specifically, the relation between the partial derivatives imposed on Assumption (D) by \cite{DeAngelis_2019_Lipschitz} excludes our gain function, but Equation \eqref{eq:H_t<_2} overcomes this issue.

\begin{proposition}[Lispchitz continuity and differentiability of the OSB]\label{pr:OSB_Lipschitz}\ \\
    The OSB $b$ is Lipschitz continuous on any closed interval of $\R_+$.
\end{proposition}

\begin{proof}
	Let $H(s, y) := W(s, y) - G(s, y)$, fix two arbitrary non-negative numbers $\underline{s}$ and $\bar{s}$ such that $\underline{s} < \bar{s}$, and consider the closed interval $I = [\,\underline{s}, \bar{s}\,]$. Proposition \ref{pr:OSB_shape} guarantees that $b$ is bounded from below and, hence, we can choose $r < \inf \lrc{b(s): s\in I}$. Then, $I\times\lrc{r}\subset\cC$, meaning that $H(s, r) > 0$ for all $s\in I$. Since $H$ is continuous (see Proposition \ref{pr:value_Lipschitz}) on $\cC$, there exists a constant $a > 0$ such that $H(s, r) \geq a$ for all $s\in I$. Therefore, for all $\delta$ such that $0 < \delta \leq a$, and all $s\in I$, there exists $y\in\R$ such that $H(s, y) = \delta$. Such a value of $y$ is unique, as $\partial_x H < 0$ on $\cC$ (see \eqref{eq:W_x}). Hence, we can denote it by $b_\delta(s)$ and define the function $b_\delta:I\rightarrow [r, b(s))$. $H$ is regular enough away from the boundary to apply the implicit function theorem, which states the differentiability of $b_\delta$ along with
    \begin{align}\label{eq:b_delta'}
        b_\delta'(s) = -\partial_t H(s, b_\delta(s)) / \partial_x H(s, b_\delta(s)).
    \end{align}
    Note that $b_\delta$ increases as $\delta \rightarrow 0$ and is upper-bounded, uniformly in $\delta$, by $b$, which is proved to be finite in Proposition \ref{pr:OSB_shape}. Hence, $b_\delta$ converges pointwise, as $\delta \rightarrow 0$, to some limit function $b_0$ such that $b_0 \leq b$ on $I$. 
    The reverse inequality follows from
    \begin{align*}
        H(s, b_0(s)) = \lim_{\delta \rightarrow 0}H(s, b_\delta(s)) = \lim_{\delta \rightarrow 0}\delta = 0,
    \end{align*}
	meaning that $(s, b_0(s))\in\cD$. Hence, $b_0 = b$ on $I$.
	
	For $(s, y)\in \cC$ such that $s\in I$ and $y > r$, consider the stopping times $\sigma^* = \sigma^*(s, y)$ and
	\begin{align*}
	    \sigma_r = \sigma_r(s, y) = \inf\{u\geq 0 : (s + u, Y_{s+u}) \notin I\times (r, \infty) \}.  
	\end{align*}
	By recalling \eqref{eq:W_t}, it readily follows that
	\begin{align}\label{eq:H_t<_1}
	    \left|\partial_t H(s, y)\right| \leq K^{(1)}\ m(s, y)
	\end{align}
	for $K^{(1)} = \max\lrc{A_1'' + 2c_0A_2' + c_0A_3'', 1}$ and
	\begin{align*}
	m(s, y) := \Esbigg{\int_0^{\sigma^*}\lrp{1 + \left|a_2''(s+u)Y_{s+u}\right|}\,\rmd u}{s, y}.
	\end{align*}
	Due to the tower property of conditional expectation, the strong Markov property, and the fact that $\sigma^*(s, y) = \sigma_r + \sigma^*\lrp{s + \sigma_r, Y_{s + \sigma_r}}$ whenever $\sigma_r \leq \sigma^*$, we have that
	\begin{align}
	    m(s, y) 
	    &= \Esbigg{\int_0^{\sigma^*\wedge\sigma_r}\lrp{1 +
	    \left|a_2''(s + u)Y_{s+u}\right|}\,\rmd u + \mathbbm{1}\lrp{\sigma_r \leq \sigma^*} m(s + \sigma_r, Y_{s + \sigma_r})}{s, y}. \label{eq:m_split}
	\end{align}
    On the set $\lrc{\sigma_r \leq \sigma^*}$, $(s + \sigma_r, Y_{s + \sigma_r}) \in \Gamma_s$ $\Pr_{s, y}$-a.s. whenever $r < y < b(s)$, with $\Gamma_s := \lrp{(s, \bar{s})\times\{r\}} \allowbreak\cup \lrp{\{\bar{s}\}\times[r, b(\bar{s})]}$. Hence, if $\sigma_r \leq \sigma^*$, then
	\begin{align}
		m\lrp{s + \sigma_r, Y_{s + \sigma_r}} &\leq \sup_{(s', y') \in \Gamma_s}m(s', y') \nonumber \\
		&\leq \sup_{(s', y') \in \Gamma_s} \Es{\int_{0}^{\infty}\lrp{1 + \left|a_2''(s' + u)Y_{s+u}\right|}\,\rmd u}{s', y'} \nonumber \\
		&\leq \sup_{(s', y') \in \Gamma_s} \lrp{\int_{0}^{\infty}\lrp{1 + \left|a_2''(s' + u)y'\right|}\,\rmd u + \int_{0}^{\infty}\Esp{\left|a_2''(s' + u)Y_u\right|}\,\rmd u} \nonumber \\
		&\leq \int_{0}^{\infty}\lrp{1 + \left|a_2''(u)M\right|}\,\rmd u + \int_{0}^{\infty}\left|a_2''(s' + u)\right|\sqrt{2u/\pi}\,\rmd u < \infty, \label{eq:m_bound}
	\end{align}
    with $M := \max\{|\sup_{s\in I}b(s)|, |r|\}$. We can guarantee the convergence of both integrals since \eqref{limit_4} implies that $\lrav{a_2''(s)}$ is asymptotically equivalent to $s^{-2}$. By plugging \eqref{eq:m_bound} into \eqref{eq:m_split}, recalling \eqref{eq:H_t<_1}, and noticing that $1 + \left|a_2''(s + u)Y_{s+u}\right| \leq 1 + A_2''M$ whenever $u\leq\sigma^*\wedge\sigma_r$, we obtain that there exists $K_I^{(2)} > 0$ such that
	\begin{align}\label{eq:H_t<_2}
	    \left|\partial_t H(s, y)\right| \leq K_I^{(2)}\ \Es{\sigma^*\wedge\sigma_r + \mathbbm{1}\lrp{\sigma_r \leq \sigma^*}}{s, y}.
	\end{align}
    Arguing as in \eqref{eq:m_split} and relying on \eqref{eq:G_y}, \eqref{eq:W_x}, and \eqref{assumption:a_2'<0}, we get that
		\begin{align}
	    |\partial_x& H(s, y)| \nonumber\\
	    &= \Es{a_2(s) - a_2(s + \sigma^*)}{s, y} = \Esbigg{\int_0^{\sigma^*}-a_2'(s + u)\,\rmd u}{s, y} \nonumber \\
	    &= \Esbigg{\int_0^{\sigma^*\wedge\sigma_r}-a_2'(s + u)\,\rmd u + \mathbbm{1}\lrp{\sigma_r \leq \sigma^*} \left|\partial_x H(s + \sigma_r, Y_{s + \sigma_r})\right|}{s, y} \nonumber \\
	    &\geq \Esbigg{\int_0^{\sigma^*\wedge\sigma_r}-a_2'(s + u)\,\rmd u + \mathbbm{1}\lrp{\sigma_r \leq \sigma^*, \sigma_r < \overline{s} - s}\left|\partial_x H(s + \sigma_r, r)\right|}{s, y}. \label{eq:H_x>_1}
	\end{align}
    Since $I\times\lrc{r}\subset\cC$, we can take $\varepsilon > 0$ such that $\cR_\varepsilon := [\,\underline{s}, \overline{s} + \varepsilon]\times(r - \varepsilon, r + \varepsilon)\subset\cC$. Thereby, $\sigma^* > \sigma_\varepsilon$ $\Pr_{s, r}$-a.s. for all $s\in I$, where
    \begin{align*}
        \sigma_\varepsilon = \sigma_\varepsilon(s, r) := \inf\lrc{u\geq 0 : \lrp{s + u, Y_{s+u}} \notin \cR_\varepsilon}.
    \end{align*}
    Hence,
	\begin{align}
	    \left|\partial_x H(s + \sigma_r, r)\right| &\geq \inf_{s\in I} \left|\partial_x H(s, r)\right| = \inf_{s\in I} \Es{a_2(s) - a_2(s + \sigma^*)}{s, r} \nonumber \\
	    &\geq \inf_{s\in I} \Es{a_2(s) - a_2(s + \sigma_\varepsilon)}{s, r} \nonumber \\
	    &\geq \inf_{s\in I} \lrp{a_2(s) - a_2(\overline{s} + \varepsilon)}\Pro{\sigma_\varepsilon = \overline{s} + \varepsilon - s}{s, r} \nonumber \\
	    &\geq \lrp{a_2(\overline{s}) - a_2(\overline{s} + \varepsilon)}\mathbb{P}\bigg(\sup_{u\leq \overline{s} + \varepsilon - \underline{s}} \lrav{Y_u} < \varepsilon\bigg) \label{eq:H_x>_2} > 0, 
	\end{align}
	where we use that $a_2$ is decreasing. Recalling that $a_2'$ is a bounded function and plugging \eqref{eq:H_x>_2} into \eqref{eq:H_x>_1}, we get that, for a constant $K_{I, \varepsilon}^{(3)} > 0$,
	\begin{align}
	    \left|\partial_x H(s, y)\right| \geq K_{I,\varepsilon}^{(3)}\ \Es{\sigma^*\wedge\sigma_r + \mathbbm{1}\lrp{\sigma_r \leq \sigma^*, \sigma_r < \overline{s} - s}}{s, y}. \label{eq:H_x>_3}
	\end{align}
    Substituting \eqref{eq:H_t<_2} and \eqref{eq:H_x>_3} into \eqref{eq:b_delta'} we get the following bound for the derivative of $b$ by some constant $K_{I,\varepsilon}^{(4)} > 0$, $y_\delta = b_\delta(s)$, and $\sigma_\delta = \sigma^*(s, y_\delta)$:
	\begin{align}
	    \left|b_\delta'(s)\right| &\leq K_{I,\varepsilon}^{(4)}\ \frac{\Es{\sigma_\delta\wedge\sigma_r + \mathbbm{1}\lrp{\sigma_r \leq \sigma_\delta}}{s, y_\delta}}{\Es{\sigma_\delta\wedge\sigma_r + \mathbbm{1}\lrp{\sigma_r \leq \sigma_\delta, \sigma_r < \overline{s} - s}}{s, y_\delta}} \nonumber \\
	    &\leq K_{I,\varepsilon}^{(4)}\lrp{1 + \frac{\Pro{\sigma_r \leq \sigma_\delta}{s, y_\delta}}{\Es{\sigma_\delta\wedge\sigma_r + \mathbbm{1}\lrp{\sigma_r \leq \sigma_\delta, \sigma_r < \overline{s} - s}}{s, y_\delta}}} \nonumber \\
	    &\leq K_{I,\varepsilon}^{(4)}\lrp{1 + \frac{\Pro{\sigma_r \leq \sigma_\delta, \sigma_r = \bar{s} - s}{s, y_\delta}}{\Es{\sigma_\delta\wedge\sigma_r}{s, y_\delta}} + \frac{\Pro{\sigma_r \leq \sigma_\delta, \sigma_r < \bar{s} - s}{s, y_\delta}}{\Es{\mathbbm{1}\lrp{\sigma_r \leq \sigma_\delta, \sigma_r < \overline{s} - s}}{s, y_\delta}}} \nonumber \\
	    &\leq K_{I,\varepsilon}^{(4)}\lrp{2 + \frac{\Pro{\sigma_r \leq \sigma_\delta, \sigma_r = \bar{s} - s}{s, y_\delta}}{\Es{\mathbbm{1}\lrp{\sigma_r\leq\sigma_\delta, \sigma_r = \bar{s} - s}\lrp{\sigma_\delta\wedge\sigma_r}}{s, y_\delta}}} \nonumber \\
	    &\leq K_{I,\varepsilon}^{(4)}\lrp{2 + \frac{1}{\bar{s} - s}}. \label{eq:|b_delta'|<}
	\end{align}
	Let $I_\varepsilon = [\,\underline{s}, \bar{s} - \varepsilon]$ for $\varepsilon > 0$ small enough. By \eqref{eq:|b_delta'|<}, there exists a constant $L_{I, \varepsilon} > 0$, independent from $\delta$, such that $|b_\delta'(s)| < L_{I, \varepsilon}$ for all $s\in I_\varepsilon$ and $0 < \delta \leq a$. Hence, Arzelà--Ascoly's theorem guarantees that $b_\delta$ converges to $b$ uniformly in $\delta \in I_\varepsilon$. 
\end{proof}

Given the local Lipschitz continuity of the OSB, proving the global continuous differentiability of the value function comes relatively easy from the law of the iterated logarithms 
and the work of \cite{DeAngelis_2020_global}, which, in turn, implies the smooth-fit condition. This approach is commented on in Remark 4.5 from \cite{DeAngelis_2019_Lipschitz}. The proposition below provides the details.   

\begin{proposition}[Global $C_1$ regularity of the value function]\label{pr:smooth-fit}\ \\
    $W$ is continuously differentiable in $\R_+\times \R$. 
\end{proposition}

\begin{proof}
    Since $W = G$ on $\cD$, and $W$ has continuous partial derivatives in $C$ (see Proposition \ref{pr:value_smoothness}), then $W$ is continuously differentiable on $\cD^\circ$ and on $\cC$, where $\cD^\circ$ stands for the interior of $\cD$. It remains to prove such regularity in $\partial \cC$ to conclude the proof.
    
    Note that the law of the iterated logarithm alongside the local Lipschitz continuity of $b$ yields the following for all $s\in\R_+$ and some constant $L_s > 0$ that depends on $s$:
    \begin{align*}
        \Pr_{s, b(s)}(\inf\{u > 0 &: Y_{s+u} > b(s+u)\} = 0) \\
        &= \lim_{\varepsilon\downarrow 0} \Pro{\inf\lrc{u > 0 : Y_{s+u} > b(s+u)} < \varepsilon}{s, b(s)} \\
        &= \lim_{\varepsilon\downarrow 0}\Pro{\sup_{u\in(0,\varepsilon)} \lrp{Y_{s+u} - b(s+u)} > 0}{s, b(s)} \\
        &= \lim_{\varepsilon\downarrow 0}\Pro{\sup_{u\in(0,\varepsilon)} \frac{Y_{s+u} - b(s+u)}{\sqrt{2u\ln(\ln(1/u))}} > 0}{s, b(s)} \\
        &\geq \lim_{\varepsilon\downarrow 0}\Pro{\sup_{u\in(0,\varepsilon)} \frac{Y_{s+u} - b(s) + L_s u}{\sqrt{2u\ln(\ln(1/u))}} > 0}{s, b(s)} \\
        &= \Pro{\limsup_{u\downarrow 0} \frac{Y_{s+u} - b(s) + L_s u}{\sqrt{2u\ln(\ln(1/u))}} > 0}{s, b(s)} = 1,
    \end{align*}
    that is, $\lrc{(s+u, Y_{s+u})}_{u\in\R_+}$ immediately enters $\cD^{\circ}$\ $\Pr_{s, b(s)}$-a.s. and, hence, Corollary 6 from \cite{DeAngelis_2020_global} guarantees that $\sigma^*(s_n, y_n) \rightarrow \sigma^*(s, b(s)) = 0\ \Pr$-a.s. for any sequence $(s_n, y_n)$ that converges to $(s, b(s))$ as $n\rightarrow\infty$.
    
    Therefore, the dominated convergence theorem and \eqref{eq:W_x} show that $$\partial_x W(s, b(s)^-) = a_2(s) = \partial_x G(s, b(s)).$$
    Since $W = G$ on $\cD$, it also holds that  
    $\partial_x W(s, b(s)^+) = \partial_x G(s, b(s)) = a_2(s)$, and, hence, $W_x$~is continuous on $\partial \cC$, which is the required smooth-fit condition.

    On the other hand, consider a sequence $s_n$ such that $(s_n, b(s))\in \cC$ for all $n$ and $s_n\uparrow s$ as $n\rightarrow\infty$. Relying again on the dominated convergence theorem and using \eqref{eq:W_t}, we get that $\partial_t W(s_n, b(s)) \rightarrow \partial_t G(s, b(s))$. We trivially reach the same convergence by taking $(s_n, b(s))\in \cD$ for all $n$, since $W = G$ on $\cD$. Arguing identically, we obtain that $\partial_t W(s_n, b(s)) \rightarrow \partial_t G(s, b(s))$ whenever $s_n \downarrow s$. Hence, $W_t$ is continuous on $\partial \cC$, which finally yields the global $C^1$ regularity of $W$.
\end{proof}

Finally, we are able to provide the solution for the OSP \eqref{eq:value_OSP}. Indeed, so far we have gathered all the regularity conditions needed to apply an extended Itô's formula to $W(s + u, Y_{s+u})$ to obtain characterizations of the value function and the OSB. The former is given in terms of an integral of the OSB, while the latter is proved to be the unique solution of a type-two, non-linear, Volterra integral equation. Both characterizations benefit from the Gaussianity of the BM, yielding relatively explicit integrands. Theorem \ref{th:OSP_solution} dives into details. Its proof needs the following lemma.

\begin{lemma}\label{lm:W->c}
    For all $(s, y)\in\R_+\times\R$,
    \begin{align*}
    \lim_{u\rightarrow\infty}\Es{W(s + u, Y_{s+u})}{s, y} = c_1 + c_0c_2,
    \end{align*}
    where $c_1$ and $c_2$ come from equations \eqref{assumption:limit_2} and \eqref{assumption:limit_1}, respectively.
\end{lemma}

\begin{proof}
    Let $s_u := s + u$ for $s, u \in\R_+$. 
    Hence, the Markov property of $Y$ implies that
    \begin{align*}
        \lim_{u\rightarrow\infty}&\,\Es{W(s_u, Y_{s_u})}{s, y} \\
        &= \lim_{u\rightarrow\infty}\Es{\sup_{\sigma\geq 0}\Es{G\lrp{s_u + \sigma, Y_{s_u + \sigma}}}{s_u, Y_{s_u}}}{s, y} \\
        &\leq \lim_{u\rightarrow\infty}\Es{\Es{\sup_{r\geq 0}G\lrp{s_u + r, Y_{s_u + r}}}{s_u, Y_{s_u}}}{s, y} \\
        &= \lim_{u\rightarrow\infty}\Es{\sup_{r\geq 0}\lrc{a_1(s_u + r) + c_0a_2(s_u + r)(s_u + r) + a_2(s_u + r)Y_{s_u+r}}}{s, y} \\ 
        &= \Esbigg{\lim_{u\rightarrow\infty}\sup_{r\geq 0}\lrc{a_1(s_u + r) + c_0a_2(s_u + r)(s_u + r) + a_2(s_u + r)Y_{s_u+r}}}{s, y} \\
        &= \Esbigg{\limsup_{u\rightarrow\infty}\lrc{a_1(s_u) + c_0a_2(s_u)s_u + a_2(s_u)Y_{s_u}}}{s, y} \\
        &= c_1 + c_0c_2,
    \end{align*}
    where the interchangeability of the limit and the mean operator is justified by the monotone convergence theorem. The last equality comes after \eqref{assumption:limit_1} and \eqref{assumption:limit_2}, along with the law of the iterated logarithm, implying that $\limsup_{u\rightarrow\infty} a_2(s_u)Y_{s_u} = 0$. 
    
    Likewise, we have that
    \begin{align*}
        \lim_{u\rightarrow\infty}\Es{W(s_u, Y_{s_u})}{s, y}
        &\geq \lim_{u\rightarrow\infty}\Es{\Es{\inf_{r\geq 0}G\lrp{s_u + r, Y_{s_u + r}}}{s_u, Y_{s_u}}}{s, y} \\
        &= \Es{\liminf_{u\rightarrow\infty}\lrc{a_1(s_u) + c_0a_2(s_u)s_u + c_0a_2(s_u)Y_{s_u}}}{s, y} \\
        &= c_1 + c_0c_2,
    \end{align*}
    which concludes the proof.
\end{proof}

\begin{theorem}[Solution of the OSP]\label{th:OSP_solution}\ \\
    The OSB related to the OSP \eqref{eq:value_OSP} satisfies the free-boundary (integral) equation
    \begin{align}\label{eq:free-boundary_eq_refined}
        G(s, b(s)) = c_1 + c_0c_2 - \int_s^\infty K(s, b(s), u, b(u)) \,\rmd u,
    \end{align}
    where the kernel $K$ is defined as
    \begin{align*}
        K(s_1, y_1, s_2, y_2) :=&\; \lrp{(a_1'(s_2) + c_0a_2(s_2) + c_0a_2'(s_2)(s_2 + y_1)}\bar{\Phi}_{s_1, y_1, s_2, y_2}\\
        &+ c_0a_2'(s_2)\sqrt{s_2 - s_1}\phi_{s_1, y_1, s_2, y_2}
    \end{align*}
    with $0 \leq s_1 \leq s_2$, $y_1, y_2\in\R$, and
    \begin{align*}
        \bar{\Phi}_{s_1, y_1, s_2, y_2} := \bar{\Phi}\lrp{\frac{y_2 - y_1}{\sqrt{s_2 - s_1}}}, \quad  
        \phi_{s_1, y_1, s_2, y_2} := \phi\lrp{\frac{y_2 - y_1}{\sqrt{s_2 - s_1}}}.
    \end{align*}
    The functions $\phi$ and $\bar{\Phi}$ are, respectively, the density and survival functions of a standard normal random variable. In addition, the integral equation \eqref{eq:free-boundary_eq_refined} admits a unique solution among the class of continuous functions $f:\R_+\rightarrow\mathbb{R}$ of bounded variation. 
    
    The value function is given by the formula
    \begin{align}\label{eq:pricing_formula_refined}
        W(s, y) &= c_1 + c_0c_2 - \int_s^\infty K(s, y, u, b(u)) \,\rmd u.
    \end{align}
\end{theorem}

\begin{proof}
    Propositions \ref{pr:value_Lipschitz}--\ref{pr:smooth-fit} provide the required regularity to apply an extended Itô's lemma (see \cite{Peskir_2005_change} for an original derivation and Lemma A2 in \cite{DAuria_2020_discounted} for a reformulation that better suits our settings) to $W(s + h, Y_{s+h})$ for $s, h\geq 0$. Since $\InfGen W = 0$ on $\cC$ and $W = G$ on $\cD$, after taking $\Pr_{s, y}$-expectation (which cancels the martingale term) it follows that
    \begin{align}
        \hspace*{-0.25cm}W(s, y) &= \Es{W(s + h, Y_{s+h})}{s, y} - \Es{\int_0^h (\InfGen W)\lrp{s + u, Y_{s+u}}\,\rmd u}{s, y} \nonumber \\
        &= \Es{W(s + h, Y_{s+h})}{s, y} - \Es{\int_0^h \partial_t G\lrp{s + u, Y_{s+u}}\Ind\lrp{Y_{s+u} \geq b(s + u)}\,\rmd u}{s, y}, \label{eq:pricing_formula_aux}
    \end{align}
    where the local-time term does not appear due to the smooth-fit condition. Hence, by taking $h\rightarrow\infty$ in \eqref{eq:pricing_formula_aux} and relying on Lemma \ref{lm:W->c}, 
    we get the following formula for the value function:
    \begin{align}
        W(s, y) &= c_1 + c_0c_2 - \Es{\int_0^\infty (\InfGen W)\lrp{s + u, Y_{s+u}}\,\rmd u}{s, y} \nonumber \\
        &= c_1 + c_0c_2 - \Es{\int_0^\infty \partial_t G\lrp{s + u, Y_{s+u}}\Ind\lrp{Y_{s+u} \geq b(s + u)}\,\rmd u}{s, y}. \label{eq:pricing_formula}
    \end{align}
    We can obtain a more tractable version of \eqref{eq:pricing_formula} by exploiting the linearity of ${y\mapsto \partial_t G(s, y)}$ (see \eqref{eq:G_t}) as well as the fact that $Y_{s+u}\sim \cN(y, u)$ under $\Pr_{s, y}$. Then,
    \begin{align*}
        \Es{Y_{s+u}\Ind\lrp{Y_{s+u}\geq x}}{s, y} = \bar{\Phi}((x - y)/\sqrt{u})y + \sqrt{u}\phi((x - y)/\sqrt{u}).
    \end{align*}
    Hence, by right-shifting the integrating variable $s$ units, we get equation \eqref{eq:pricing_formula_refined}.
    
    Take now $y\downarrow b(s)$ in both \eqref{eq:pricing_formula} and \eqref{eq:pricing_formula_refined} to derive the free-boundary equation
    \begin{align}\label{eq:free-boundary_eq}
        G(s, b(s)) &= c_1 + c_0c_2 - \Es{\int_0^\infty \partial_t G\lrp{s + u, Y_{s+u}}\Ind\lrp{Y_{s+u} \geq b(s + u)}\,\rmd u}{s, b(s)},
    \end{align}
    alongside its more explicit expression \eqref{eq:free-boundary_eq_refined}.
    
    The uniqueness of the solution of equation \eqref{eq:free-boundary_eq} follows a well-known methodology first developed by \cite[Theorem 3.1]{Peskir_2005_American} that we omit here for the sake of briefness.
\end{proof}

\section{Solution of the original OSP}\label{sec:solutionoriginal}

In this section we continue with the notation used in Section \ref{sec:formulation}.

Recall that Proposition \ref{pr:OSP_equivalence} dictates the equivalence between the OSPs \eqref{eq:value_OSP_GMB} and \eqref{eq:value_OSP_BM}, and gives explicit formulae to link their value functions and OSTs. Consequently, it follows that the stopping time $\tau^*(t, x)$ defined in Proposition \ref{pr:OSP_equivalence} in terms of $\sigma^*(s, y)$ is not only optimal for \eqref{eq:value_OSP_GMB}, but it holds the following representation under $\PPr_{t, x}$:
\begin{align}\label{eq:OST_OSB}
    \tau^*(t, x)
    &= \inf\lrc{u\geq 0 : X_{t + u} \geq \rmb_{T, z}(t + u)},\quad \rmb_{T, z}(t) := G_{T, z}(s, b_{T, z}(s)),
\end{align}
where $\rmb_{T, z}$ and $b_{T, z}$ are, respectively, the OSBs related to \eqref{eq:value_OSP_GMB} and \eqref{eq:value_OSP_BM}, and $s$ is defined, in terms of $t$, in Proposition \ref{pr:OSP_equivalence}. Note that $b_{T, z}$ coincides with the function defined in \eqref{eq:OSB}, with constants $c_0$, $c_1$, and $c_2$, from \eqref{eq:gain_function}, \eqref{assumption:limit_1}, and \eqref{assumption:limit_2}, taking the values
\begin{align}\label{eq:constants}
        c_0 = z - \alpha(T),\quad c_1 = \alpha(T),\quad c_2 = 1,
\end{align}
where $\alpha$ comes from \eqref{pr:GMB_to_BM_3}, Proposition \ref{pr:GMB} (see also Remark \ref{rm:c0-a1-a2}).

Moreover, it is not necessary to compute $W_{T, z}$ and $b_{T, z}$ to obtain $V_{T, z}$ and $\rmb_{T, z}$. By considering the infinitesimal generator of $\lrc{\lrp{t, X_t}}_{t \in [0, T]}$, $\rmL$, letting $s_\varepsilon = s + \varepsilon$ and $t_\varepsilon = \gamma_T^{-1}(s_\varepsilon)$ for $\varepsilon > 0$, and using \eqref{eq:equivalence_value} alongside the chain rule, we get that
\begin{align}
    \lrp{\InfGen W_{T, z}}(s, y) :=&\; \lim_{\varepsilon\rightarrow 0}\varepsilon^{-1}\lrp{\Es{W_{T, z}\lrp{s_\varepsilon, Y_{s_\varepsilon}}}{s, y} - W_{T, z}(s, y)} \nonumber \\
    =&\; \lim_{\varepsilon\rightarrow 0}\varepsilon^{-1}\lrp{\EEs{V_{T, z}(t_\varepsilon, X_{t_\varepsilon})}{t, x} - V_{T, z}(t, x)} \nonumber \\
    =&\; \lrp{\rmL V_{T, z}}(t, x)\big[\gamma_T^{-1}\big]'(s). \label{eq:infinitesimal_relation}
\end{align}
We recall the relation between $s$ and $t$, and $y$ and $x$, in Proposition \ref{pr:OSP_equivalence}. After integrating with respect to $\gamma_T^{-1}(u)$ instead of $u$ in \eqref{eq:pricing_formula_aux}, keeping in mind \eqref{eq:constants} and \eqref{eq:infinitesimal_relation}, and recalling that $\rmL V_{T, z}(t, x) = 0$ for all $x\leq \rmb_{T, z}(t)$ and $V_{T, z}(t, x) = x$ for all $x\geq\rmb_{T, z}(t)$, we get the formula
\begin{align}
    V_{T, z}(t, x) &= z - \EEs{\int_0^{T-t} (\rmL V_{T, z})(t + u, X_{t + u})\,\rmd u}{t, x} \nonumber \\
    &= z - \EEs{\int_0^{T-t} \mu(t + u, X_{t + u})\mathbbm{1}(X_{t + u} \geq \rmb_{T, z}(t + u))\,\rmd u}{t, x}, \label{eq:pricing_formula_original}
\end{align}
where, in alignment to \eqref{eq:tdOUB_coeff_BM},
\begin{align*}
    \mu(t, x) :=&\; \lim_{u\downarrow 0}u^{-1}\EEs{X_{t + u} - x}{t, x} 
    = \theta(t)(\kappa(t) - x)\\
    =&\; \alpha'(t) + \lrp{x - \alpha(t)}\frac{\beta_T'(t)}{\beta_T(t)} + (z - \alpha(T))\beta_T(t)\gamma_T'(t).
\end{align*}
As we did to obtain \eqref{eq:pricing_formula_refined}, the linearity of $x\mapsto\mu(t, x)$ and the Gaussian marginal distributions of $X$, allow us to produce a refined version of \eqref{eq:pricing_formula_original}:
\begin{align}\label{eq:pricing_formula_original_refined}
    V_{T, z}(t, x) = z - \int_t^T\rK(t, x, u, \rmb_{T, z}(u))\,\rmd u,
\end{align}
where
\begin{align}
    \rK(t_1,& x_1, t_2, x_2) \nonumber \\
    :=&\; \theta(t_2)\lrp{(\kappa(t_2) - \EEs{X_{t_2}}{t_1, x_1})\rPhi_{t_1, x_1, t_2, x_2} - \sqrt{\VVs{X_{t_2}}{t_1}}\frac{\beta_T'(t_2)}{\beta_T(t_2)}\rphi_{t_1, x_1, t_2, x_2}} \label{eq:kernel_sde} \\
    =&\; \lrp{\alpha'(t_2) + \lrp{\EEs{X_{t_2}}{t_1, x_1} -  \alpha(t_2)}\frac{\beta_T'(t_2)}{\beta_T(t_2)} + (z - \alpha(T))\beta_T(t_2)\gamma_T'(t_2)}\rPhi_{t_1, x_1, t_2, x_2} \nonumber\\
    & + \sqrt{\VVs{X_{t_2}}{t_1}}\frac{\beta_T'(t_2)}{\beta_T(t_2)}\rphi_{t_1, x_1, t_2, x_2}, \label{eq:kernel_BM}
\end{align}
with $0 \leq t_1 \leq t_2 < T$, $x_1, x_2\in\R$, and
\begin{align*}
    \rPhi_{t_1, x_1, t_2, x_2} := \bar{\Phi}\lrp{\frac{x_2 - \EEs{X_{t_2}}{t_1, x_1}}{\sqrt{\VVs{X_{t_2}}{t_1}}}}, \quad  
    \rphi_{t_1, x_1, t_2, x_2} := \phi\lrp{\frac{x_2 - \EEs{X_{t_2}}{t_1, x_1}}{\sqrt{\VVs{X_{t_2}}{t_1}}}},
\end{align*}
and, as stated in \eqref{eq:tdOU_mean_1}, \eqref{eq:tdOU_cov_fac}, and \eqref{eq:tdOUB_coeff_BM},
\begin{align}
    \EEs{X_{t_2}}{t_1, x_1} &= \varphi(t_2)\lrp{\frac{x_1}{\varphi(t_1)} + \int_{t_1}^{t^2}\frac{\kappa(u)\theta(u)}{\varphi(u)}\,\rmd u} \label{eq:mean_sde}\\
    &= \alpha(t_2) +\beta_T(t_2)\lrp{(z - \alpha(T))\gamma_T(t_2) - \frac{x_1 - \alpha(t_1) - \beta_T(t_1)\gamma_T(t_1)(z - \alpha(T))}{\beta_T(t_1)}}, \nonumber \\
    \VVs{X_{t_2}}{t_1} &= \varphi^2(t_2)\int_{t_1}^{t_2}\frac{\nu^2(u)}{\varphi^2(u)}\,\rmd u \label{eq:var_sde} \\
    &= \beta_T(t_1)\gamma_T(t_1)\beta_T(t_2), \nonumber
\end{align}
with $\varphi(t) = \exp\lrc{-\int_{0}^{t}\theta(u)\,\rmd u}$. Hence, after taking $x\downarrow\rmb(t)$ in \eqref{eq:pricing_formula_original} (or by directly expressing \eqref{eq:free-boundary_eq} in terms of the original OSP, as we did to obtain \eqref{eq:pricing_formula_original} from \eqref{eq:pricing_formula}), we get the free-boundary equation
\begin{align*}
    \rmb_{T, z}(t) &= z - \EEs{\int_0^{T-t} (\InfGen_X V_{T, z})(t + u, X_{t + u})\,\rmd u}{t, \rmb_{T, z}(t)} \\
    &= z - \EEs{\int_0^{T-t} \mu(t + u, X_{t + u})\mathbbm{1}(X_{t + u} \geq \rmb_{T, z}(t + u))\,\rmd u}{t, \rmb_{T, z}(t)},
\end{align*}
which is also expressible as
\begin{align}\label{eq:free-boundary_eq_original_refined}
    \rmb_{T, z}(t) = z - \int_t^T \rK(t, \rmb_{T, z}(t), u, \rmb_{T, z}(u))\,\rmd u.
\end{align}
The uniqueness of the solution of the Volterra-type integral equation \eqref{eq:free-boundary_eq_original_refined} follows from that of \eqref{eq:free-boundary_eq_refined}.

\begin{remark}
    We highlight some smoothness properties that the value function~$V$ and the OSB $\rmb$ inherit from $W$ and $b$, based on the equivalences \eqref{eq:equivalence_value} and \eqref{eq:OST_OSB}. 
    
    From the Lipschitz continuity of $W$ on compact sets of $\R_+\times\R$ (see Proposition~\ref{pr:value_Lipschitz}), it follows that of $V$ in compact sets of $[0, T)\times\R$. Higher smoothness of $V$ is also attained away from the boundary, $(t, \rmb(t))$ for all $t\in[0, T)$, as it follows from Proposition \ref{pr:value_smoothness}. The continuous differentiability of $W$ obtained in Proposition \ref{pr:smooth-fit} implies that of $V$.
    
    The OSB $\rmb$ is Lipschitz continuous on any closed subinterval of $[0, T)$, which is a consequence of Proposition \ref{pr:OSB_Lipschitz}.
\end{remark}

\section{Numerical results}\label{sec:numerical_results}

In this section we shed light on the OSB's shape by using a Picard iteration algorithm to solve the free-boundary equation \eqref{eq:free-boundary_eq_original_refined}. This approach is commonly used in the optimal stopping literature; see, e.g., the works of \cite{Detemple_2020_value} and \cite{DeAngelis_2020_optimal}.

A Picard iteration scheme approaches \eqref{eq:free-boundary_eq_original_refined} as a fixed-point problem. From an initial candidate boundary, it produces a sequence of functions by iteratively computing the integral operator in the right-hand side of \eqref{eq:free-boundary_eq_original_refined}, until the error between consecutive boundaries is below a prescribed threshold. More precisely, for a partition $0 = t_0 < t_1 < \cdots < t_N = T$ of $[0, T]$, $N\in\mathbb{N}$, the updating mechanism that generates subsequent boundaries follows after the discretization of the integral in \eqref{eq:free-boundary_eq_original_refined} by using a right Riemann sum:
\begin{align}
    \rmb_i^{(k)} &= z - \sum_{j = i}^{N - 2} \rK\lrp{t_i, \rmb_i^{(k - 1)}, t_{j + 1}, \rmb_{j + 1}^{(k - 1)}}(t_{j + 1} - t_j), \quad i = 0, 1, \dots, N-2, \label{eq:picard1} \\
    \rmb_{N-1}^{(k)} &= \rmb_N^{(k)} = z, \label{eq:picard2}
\end{align}
for $k = 1, 2, \dots$ and with $\rmb_i^{(k)}$ standing for the value of the boundary at $t_i$ output after the $k$-th iteration. We neglect the $(N-1)$-addend of the sum, and instead consider \eqref{eq:picard2}, since $\rK(t, x, T, z)$ is not well defined. As the integral in \eqref{eq:free-boundary_eq_original_refined} is finite, the last piece vanishes as $t_{N - 1}$ approaches $T$. Given that $\rmb(T) = z$, we set the initial constant boundary $\rmb_i^{(0)} = z$ for all $i = 0, \dots, N$. We stop the fixed-point algorithm when the relative (squared) $L_2$-distance between the consecutive discretized boundaries, defined~as 
\begin{align*}
    d_k := \frac{\sum_{i=1}^N \lrp{\rmb_i^{(k)} - \rmb_i^{(k-1)}}^2(t_i - t_{i-1})}{\sum_{i=1}^N \lrp{\rmb_i^{(k)}}^2(t_i - t_{i-1})},
\end{align*} 
is lower than $10^{-3}$. 

We show empirical evidence of the convergence of this Picard iteration scheme in Figures \ref{fig:algorithm_validation}--\ref{fig:coeff_effect}. For each computer drawing of the OSB, we provide smaller images at the bottom with the (logarithmically-scaled) errors $d_k$, which tend to decrease at a steep pace, making the algorithm converge ($d_k < 10^{-3}$) after few iterations.

We perform all boundary computations by relying on the SDE representation of the kernel $\rK$ defined at \eqref{eq:kernel_sde}, \eqref{eq:mean_sde}, and \eqref{eq:var_sde}, since we adopted the viewpoint of a GMB derived from conditioning a time-dependent OU process to degenerate at the horizon. The relation between the ``parent'' OU process and the resulting OUB is neatly stated in \cite[Section 3]{Buonocore_2013_some}, although we include here a modified version that fits our notation better. That is, if $\wt{X} = \{\wt{X}_t\}_{t\in[0, T]}$ solves the SDE
\begin{align}\label{eq:OU_SDE}
    \rmd \wt{X}_t = \wt{\theta}(t)(\wt{\kappa}(t) - \wt{X}_t)\,\rmd t + \wt{\nu}(t)\,\rmd B_t,\quad t\in[0, T],
\end{align}
then, the corresponding GMB is an OUB that solves the SDE
\begin{align}\label{eq:OUB_SDE}
    \rmd X_t = \theta(t)(\kappa(t) - X_t)\,\rmd t + \nu(t)\,\rmd B_t,\quad t\in(0, T),
\end{align}
with
\begin{align}\label{eq:OU_OUB_coeff}
    \left\{
    \begin{aligned}
        \theta(t) &= \wt{\theta}(t) + \frac{\wt{\nu}^2(t)}{\wt{\varphi}^2(t)\int_t^T \wt{\nu}^2(u)/\wt{\varphi}(u)\,\rmd u}, \\
        \kappa(t) &= \wt{\kappa}(t) + \frac{\wt{\nu}^2(t)}{\theta(t)}\frac{x - \wt{\varphi}(T)\int_t^T \wt{\kappa}(u)\wt{\theta}(u)/\wt{\varphi}(u)\,\rmd u}{\wt{\varphi}(t)\wt{\varphi}(T) \int_t^T \wt{\nu}^2(u)/\wt{\varphi}(u)\,\rmd u}, \\
        \nu(t) &= \wt{\nu}(t),
    \end{aligned}
    \right.
\end{align}
and where $\wt{\varphi}(t) = \exp\{-\int_0^t\wt{\theta}(u)\,\rmd u\}$. We choose representations \eqref{eq:OU_SDE} and \eqref{eq:OUB_SDE} for GM processes and GMBs, over those given in Lemma \ref{lm:GM_def} and \eqref{def:GMB_BM} from Proposition \ref{pr:GMB}, as they have a more intuitive meaning. Indeed, recall that $\theta$ ($\wt{\theta}$) indicates the strength with which the underlying process is pulled towards the mean-reverting level $\kappa$ ($\wt{\kappa}$), while $\nu$ ($\wt{\nu}$) regulates the intensity of the white-noise.

Figure \ref{fig:algorithm_validation} shows the numerically computed OSB when the underlying diffusion is a BB, that is, when $\wt{\theta}(t) = 0$ and $\wt{\nu}(t) = \sigma$, for all $t\in[0, T]$ and $\sigma > 0$. We rely on such a case to empirically validate the Picard algorithm's accuracy in Figure \ref{fig:algorithm_validation}(a) by comparing against the explicit OSB of a BB, which is known to take the form $z + K\sigma\sqrt{T - t}$, for $K\approx 0.8399$. 
This result was originally due to \cite{Shepp_1969_explicit}. Notice in Figure \ref{fig:algorithm_validation}(b) how the numerical boundary approaches the real one as the time partition becomes thinner.

For all boundary computations, $T = 1$ and $N = 500$ were set unless otherwise stated. We used the logarithmically-spaced partition $t_i = \ln\lrp{1 + i(e - 1)/N}$, since numerical tests suggested that the best performance is achieved when using a non-uniform mesh whose distances $t_i - t_{i-1}$ smoothly decrease. Figure \ref{fig:algorithm_validation}(c) illustrates such an effect of the mesh increments by comparing the performance of the logarithmically-spaced partition against an equally-spaced one and another that is also equally spaced until the second last node, where it suddenly shrinks the distance to a fourth of the regular space. Note how the first partition significantly outperforms the other two with a lower overall $L_2$-error due to its better accuracy near the horizon. Intuition might dictate that introducing the sudden shrink at the horizon may result in better performance by diminishing the error that arises when considering \eqref{eq:picard2}, yet Figure \ref{fig:algorithm_validation}(c) indicates otherwise.

\begin{figure}[ht]
    \centering
    \subfloat[$\wt{\theta} \equiv 0$.]{%
    \resizebox*{0.32\textwidth}{!}{\includegraphics{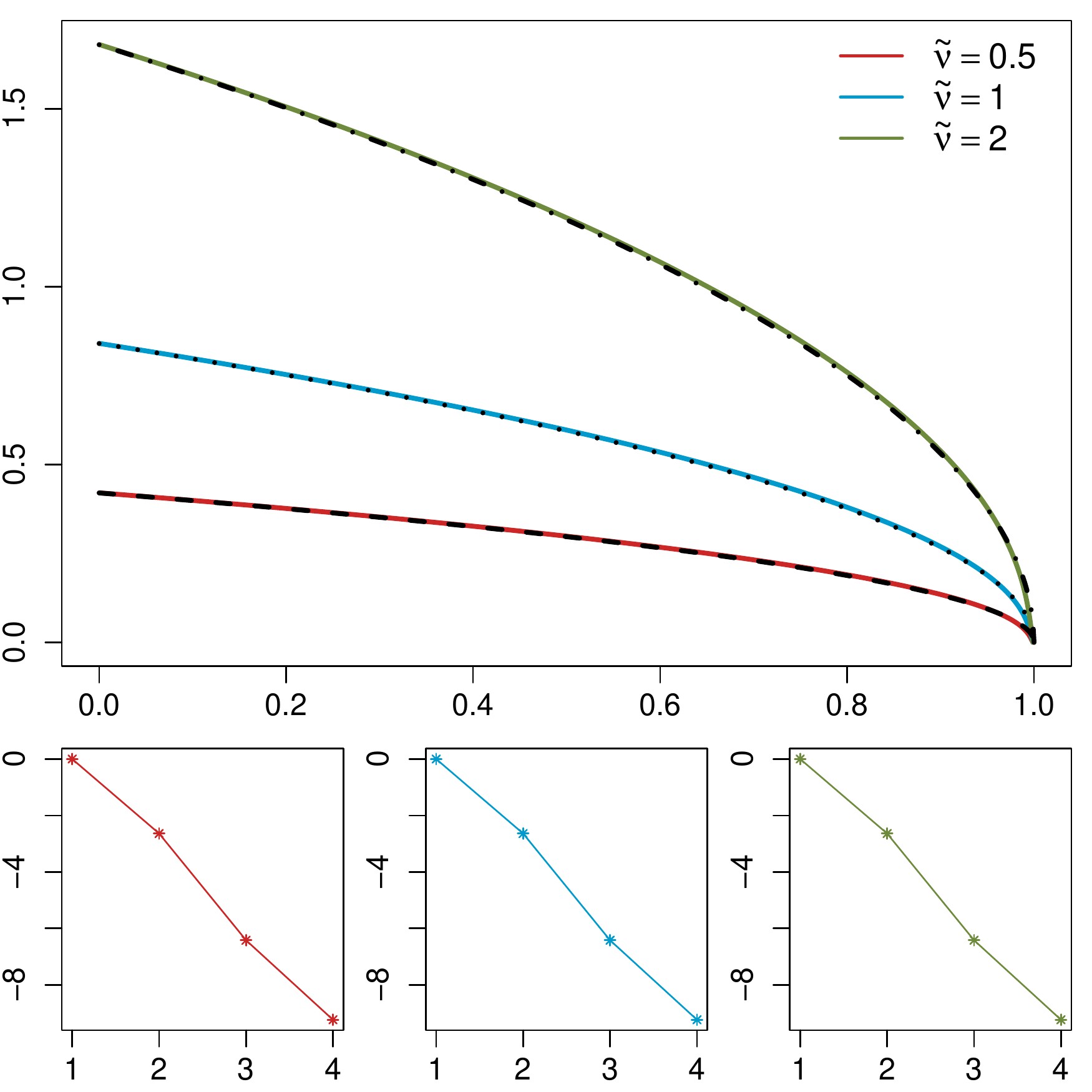}}}\hspace{2pt}
    \subfloat[$\wt{\theta} \equiv 0$, $\wt{\nu} \equiv 1$.]{%
    \resizebox*{0.32\textwidth}{!}{\includegraphics{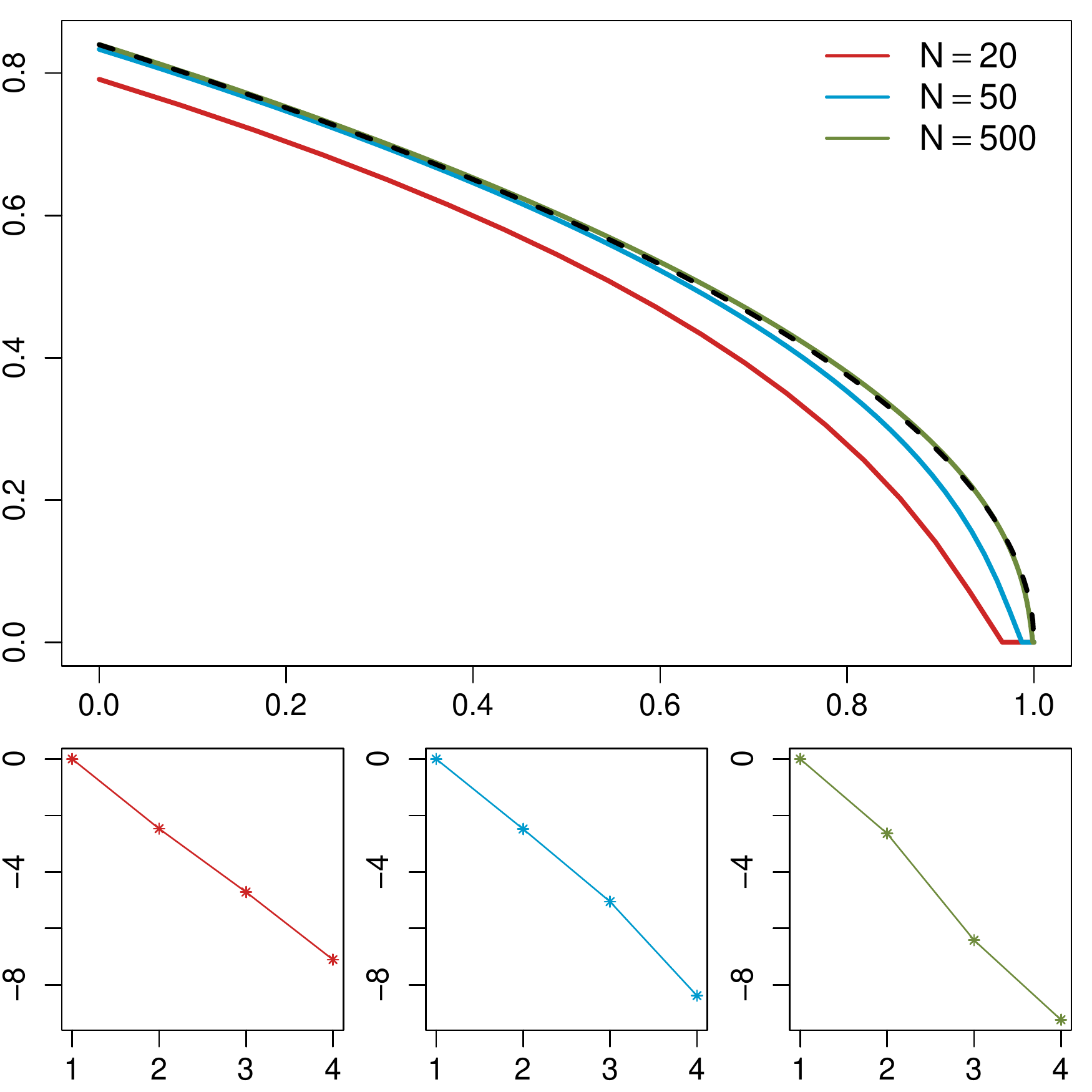}}}\hspace{2pt}
    \subfloat[$\wt{\theta} \equiv 0$, $\wt{\nu} \equiv 1$, $N = 20$.]{%
    \resizebox*{0.32\textwidth}{!}{\includegraphics{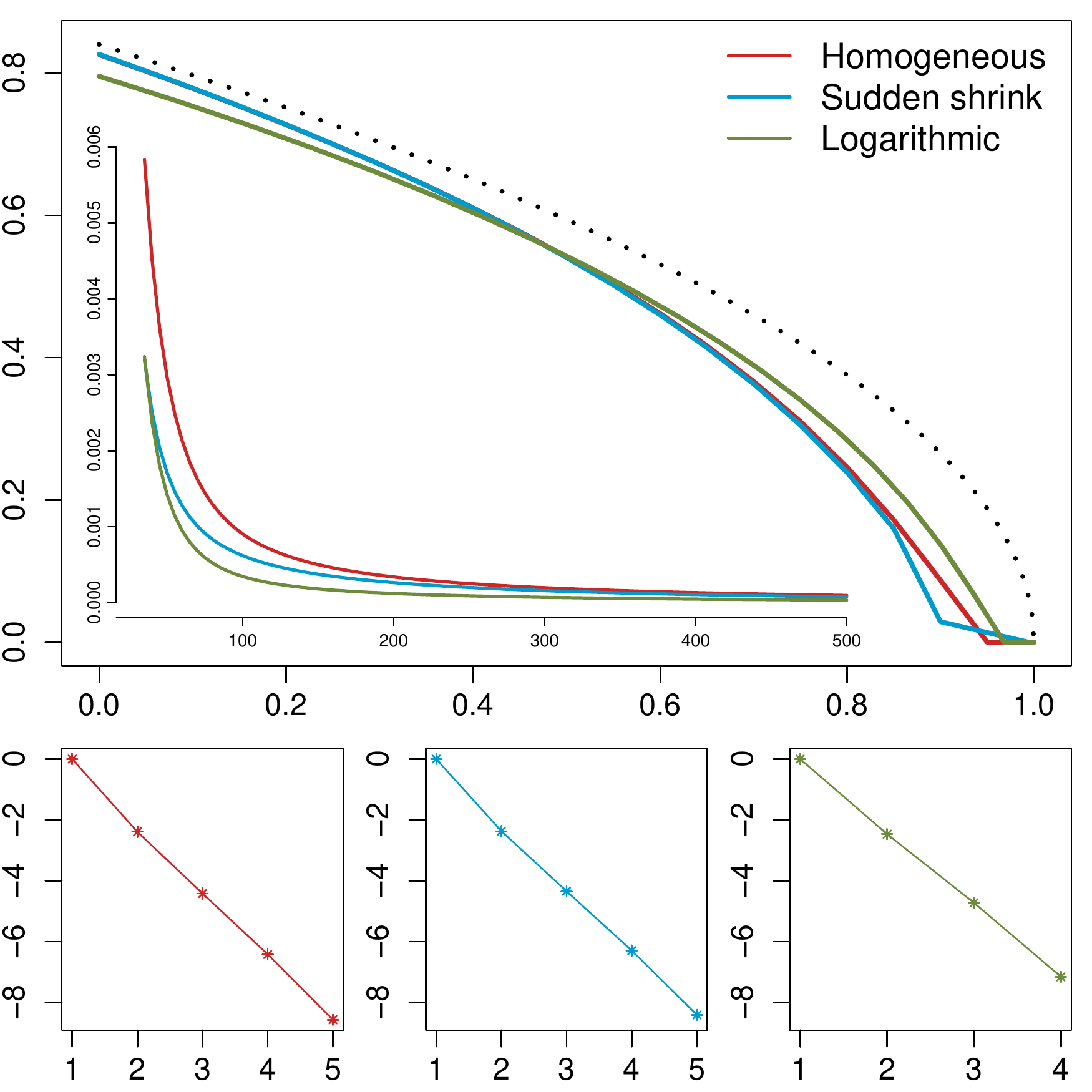}}}
    \caption{\small The picture shows a comparison between the exact OSB of a BB and its numerical computation, which is obtained by setting $\wt{\theta}\equiv 0$ and taking a constant volatility $\wt{\nu}$ in the OU representation \eqref{eq:OU_SDE}. For the images on top, the solid colored lines represent the computed OSBs for the different choices of the volatility coefficient $\wt{\nu}$ (image (a)), the partition length $N$ (image (b)), and the type of partition considered (image (c)). Black dashed, dotted, and dashed-dotted lines stand for the OSB of a BB associated with the different values of $\wt{\nu}$. Specifications are shown in the legend and caption of each image. Image (c) accounts for a subplot that shows, as a function of the partition size $N$ ($x$ axis), the evolution of the relative $L_2$ error between the different computed boundaries and the true one ($y$ axis). The smaller images at the bottom display the log-errors $\log_{10}(d_k)$ between consecutive boundaries for each iteration $k = 1, 2,\ldots$ of the Picard algorithm.} 
    \label{fig:algorithm_validation}
\end{figure}

Figure \ref{fig:coeff_effect} shows the numerical computation of OSBs for more general cases rather than the BB. It shows how changing the coefficients of the process affects the OSB shape. In the first two rows of images, we visually represent the transformation of coefficients~\eqref{eq:OU_OUB_coeff}. The volatility is excluded as it remains the same after ``bridging'' the OU process. To compare the slopes we rely on $1/\wt{\theta}(t)$ and $1/\theta(t)$, as $\theta(t) \rightarrow \infty$ as $t\rightarrow \infty$ (see \eqref{def:GMB_OUB} in Proposition \ref{pr:GMB}) and, thus, its explosion would have obscured the shape of the bounded function $\wt{\theta}$, had they been plotted in the same graph. In alignment with the meaning behind each time-dependent coefficient, the OSB is pulled towards $\wt{\kappa}$ with a strength directly proportional to $\wt{\theta}$. This pulling force conflicts with the much stronger one towards the pinning point of the bridge process, resulting in an attraction towards the ``bridged'' mean-reverting level $\kappa$ with strength dictated by $\theta$. We recall that modifying $\wt{\nu}$, and thus $\nu$, is equivalent to change $\theta$ due to \eqref{eq:theta_sigma}. We remind that the functions $\Phi$ and $\phi$ in Figure \ref{fig:coeff_effect} stand for the distribution and the density of a standard normal random variable. The former is used to smoothly represent sudden changes of regimes, while the latter introduces smooth temporal anomalies. For instance, $\wt{\kappa}(t) = 2\Phi(50t - 25) - 1$ rapidly changes the mean-reverting level of the underlying process from $-1$ to $1$ around $t = 0.5$, and $\wt{\nu}(t) = 1 + \sqrt{2\pi}\phi(100t - 25)$ introduces a brief period of increased volatility around $t = 0.25$, before and after which the volatility remains at (constant) baseline levels. Periodic fluctuations of the parameters were also considered, as they typically arise in problems that account for seasonality.

\begin{figure}[ht]
	\centering
	\begin{subfigure}[b]{0.32\textwidth}
		\includegraphics[width = \textwidth, height = 8cm, keepaspectratio]{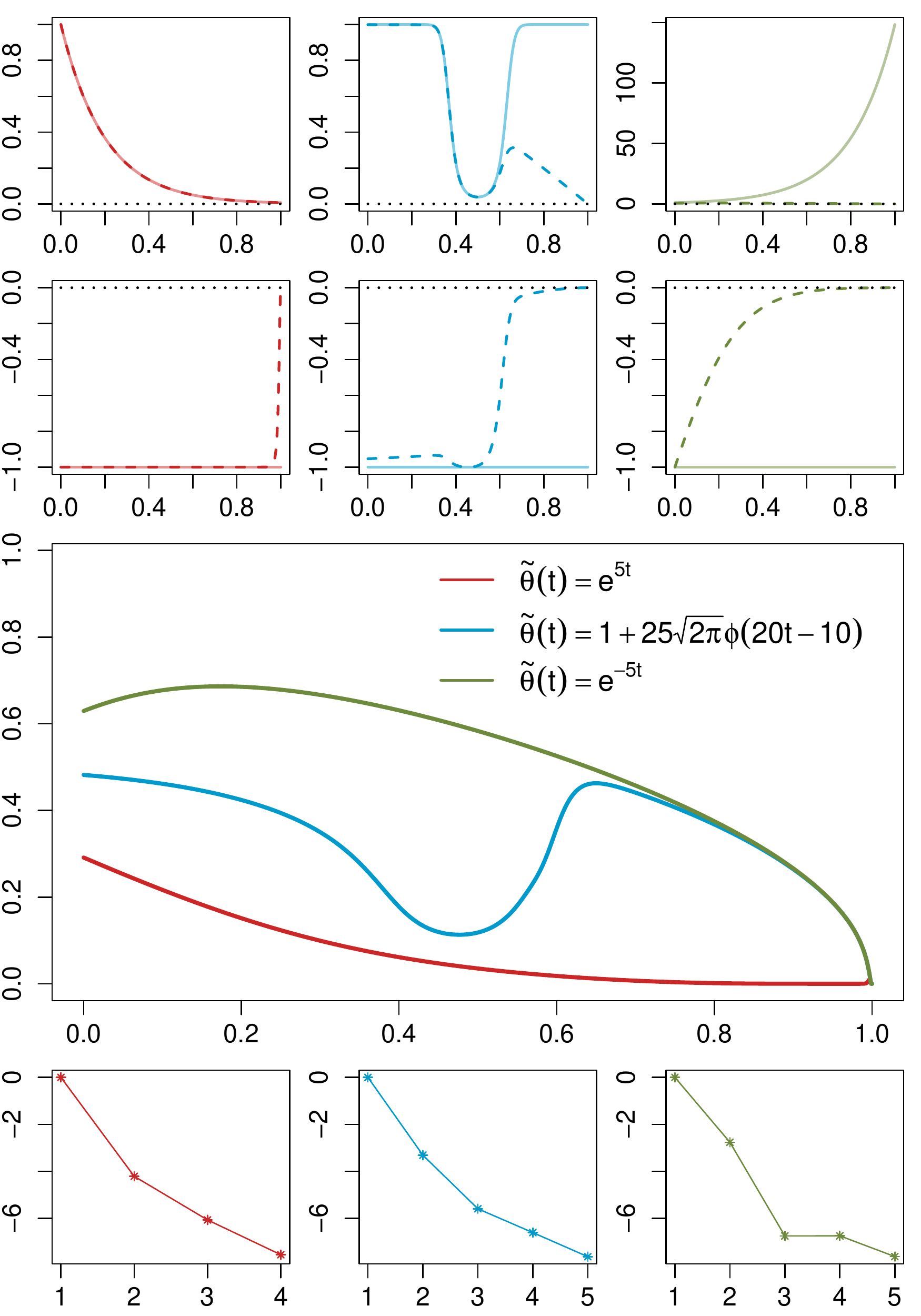}
		\subcaption{$\wt{\kappa} \equiv -1$, $\wt{\nu} \equiv 1$.}
	\end{subfigure}
	\begin{subfigure}[b]{0.32\textwidth}
		\includegraphics[width = \textwidth]{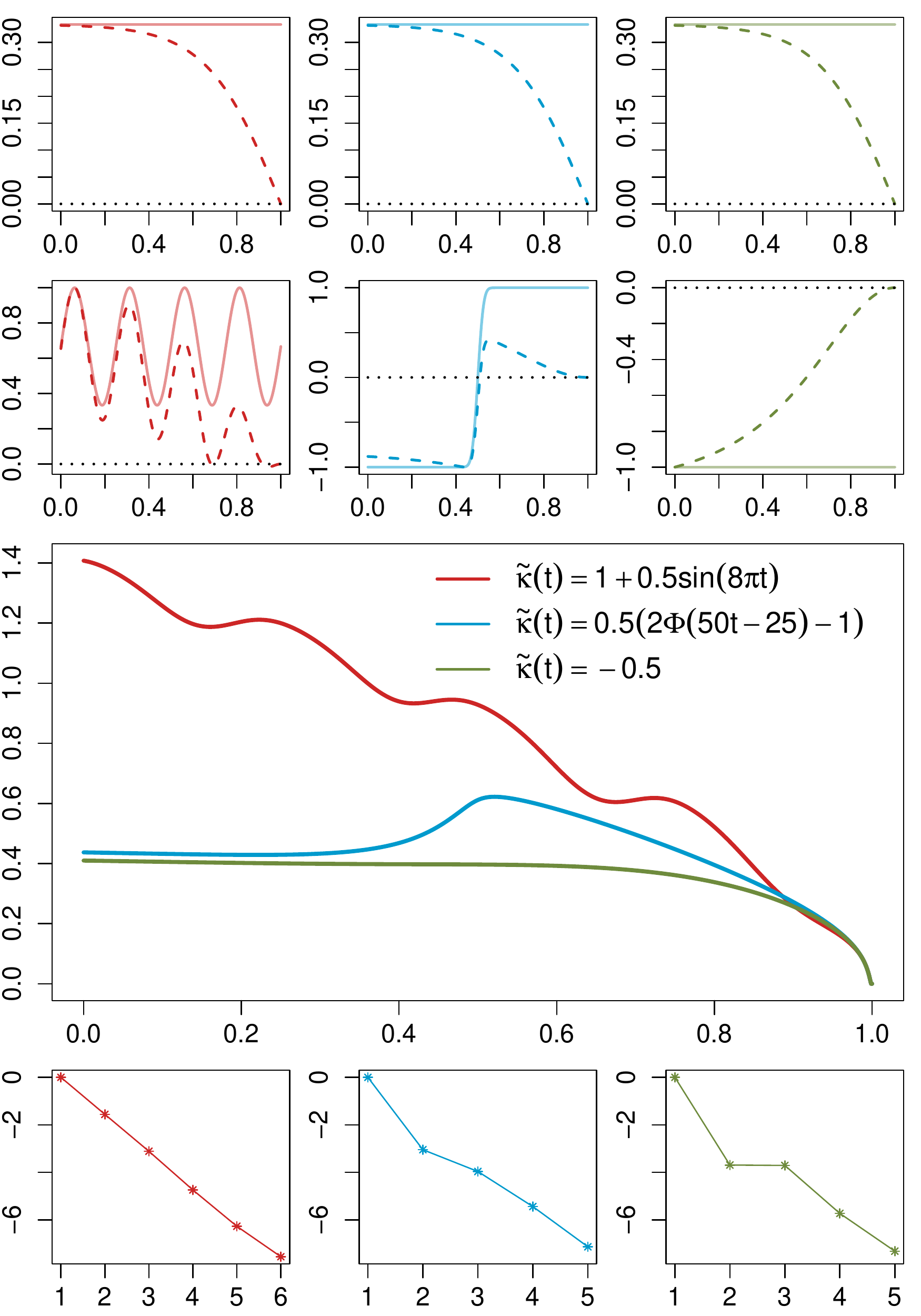}
		\subcaption{$\wt{\theta} \equiv 3$, $\wt{\nu} \equiv 1$.}
	\end{subfigure}
	\begin{subfigure}[b]{0.32\textwidth}
		\includegraphics[width = \textwidth]{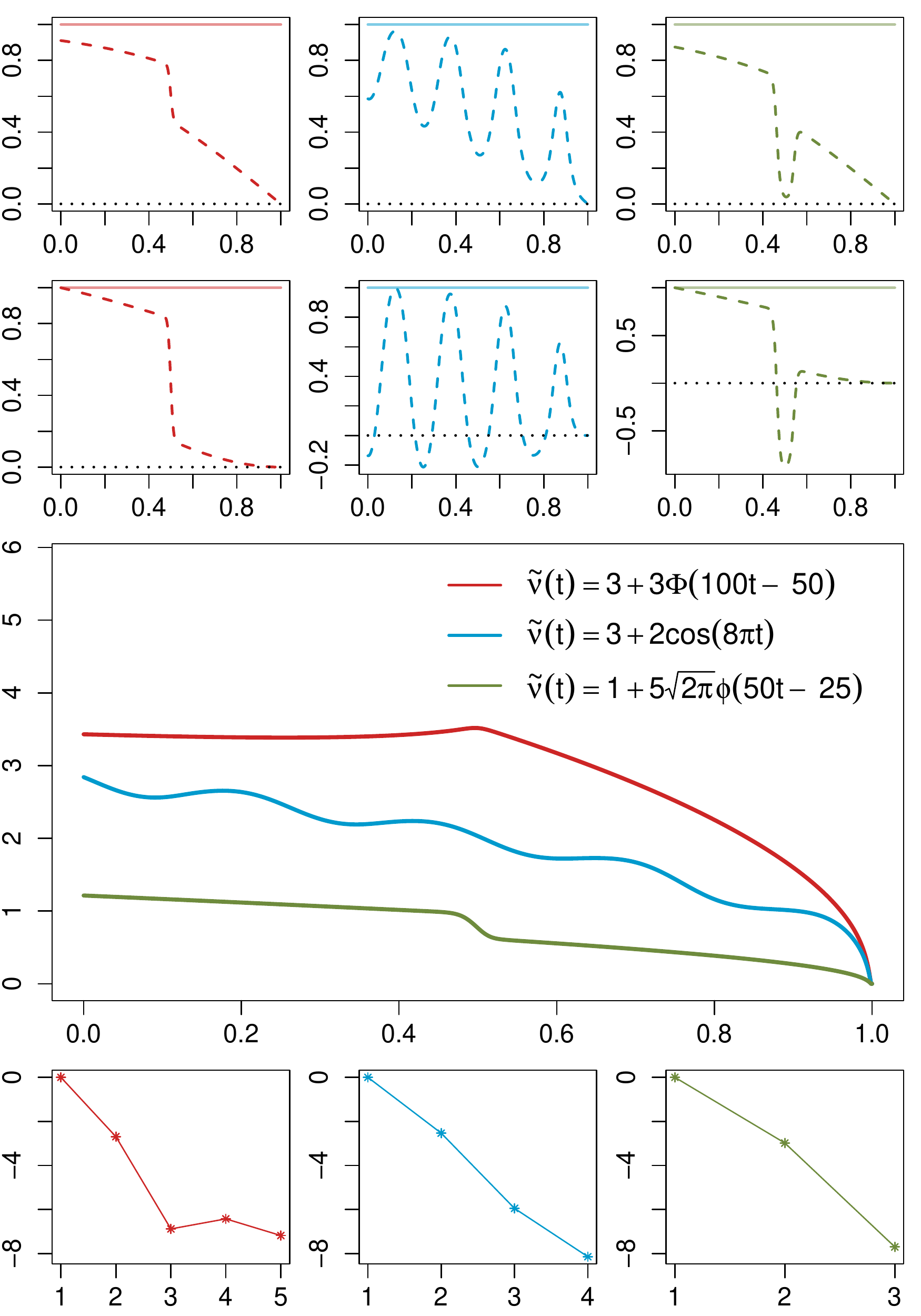}
		\subcaption{$\wt{\theta} \equiv 1$, $\wt{\kappa} \equiv 1$.}
	\end{subfigure}
	\caption{\small The first row of three plots shows $1/\wt{\theta}$ (continuous line) versus $1/\theta$ (dashed line) for the different choices of the slope $\wt{\theta}$ (image (a)), the mean-reverting level $\wt{\kappa}$ (image (b)), and the volatility $\wt{\nu}$ (image~(c)) functions. Specifications of the functions are given in the legend and caption of each image. The second row does the same for $\wt{\kappa}$ and $\kappa$. The main plot, in the third row, shows in solid colored lines the computed OSBs. The smaller images at the bottom display the log-errors $\log_{10}(d_k)$ between consecutive boundaries for each iteration $k = 1, 2,\ldots$ of the Picard algorithm.}
	\label{fig:coeff_effect}
\end{figure}


Notice that, after Proposition \ref{pr:GMB}, it readily follows that all coefficients $\theta$, $\kappa$, and $\nu$ used in this section meet assumptions \eqref{assumption:continuous_diff}--\eqref{assumption:a_2'<0}, as they are twice continuous differentiable, $\theta(t) > 0$ for all $t\in[0, T)$, and satisfy conditions \eqref{eq:theta_explotion} and \eqref{eq:theta_sigma}.

The R code in the public repository \url{https://github.com/aguazz/OSP_GMB} implements the Picard iteration algorithm \eqref{eq:picard1}--\eqref{eq:picard2}. The repository allows for full replicability of the above numerical examples. 

\section{Concluding remarks}\label{sec:conclusions}

We solved the finite-horizon OSP of a GMB by proving that its OSB uniquely solves the Volterra-type integral equation \eqref{eq:free-boundary_eq_original_refined}.  

GMBs were comprehensively studied in Section \ref{sec:GMB}, where four equivalent definitions were presented, making it easier to identify, create, and understand them from different perspectives. One of these representations allows bypassing the challenge of working with diffusions with non-bounded drifts and, instead, working with an equivalent infinite-horizon OSP with a BM underneath. Equations \eqref{eq:OST_OSB} explicitly relate both OSTs and OSBs, while \eqref{eq:pricing_formula_original_refined} and \eqref{eq:free-boundary_eq_original_refined} give the value formula and free-boundary equation in the original~OSP.

The method for solving the alternative OSP consisted in solving the associated free-boundary problem. To do so, several regularity properties about the value function and the OSB were obtained in Section \ref{sec:solution}, among which the local Lipschitz continuity of the OSB stands out as a remarkable property.

We approached the free-boundary equation as a fixed-point problem in Section \ref{sec:numerical_results} to numerically explore the geometry of the OSB. This provided insights about its shape for different sets of coefficients of the underlying GMB, seen as bridges derived from conditioning a time-dependent OU process to hit a pinning point at the horizon. The OSB shows an attraction toward the mean-reverting level, which fades away as time approaches the horizon, where the boundary hits the OUB's pinning~point.

In the context of gain functions beyond the identity, it is worth noting that the representation \eqref{eq:GMB_to_BM} can still be used to transform the initial OSP into an infinite-horizon one with a BM underneath. This prompts the question of extending the methodology in Section \ref{sec:solution} to address more flexible gain functions. A practical starting point for this extension might be considering a space-linear gain function, which results in simple forms for its partial derivatives (recall \eqref{eq:G_t} and \eqref{eq:G_y}) and keeps available the comparison method used in Proposition \ref{pr:OSB_shape} to obtain the boundedness of the OSB. Also, the new gain function should account for boundedness and time-wise differentiability regularities equivalent to Assumptions \eqref{assumption:continuous_diff}--\eqref{assumption:a_2'<0}.

\section*{Acknowledgements}

The authors thank the anonymous referees for their comments, which helped in improving the quality of the manuscript.

\section*{Funding}

The authors acknowledge support from grants PID2020-116694GB-I00 (first and second authors), and PID2021-124051NB-I00 (third author), funded by MCIN/AEI/10.13039/\-501100011033 and by ``ERDF A way of making Europe''.
The third author acknowledges support from ``Convocatoria de la Universidad Carlos III de Madrid de Ayudas para la recualificación del sistema universitario español para 2021--2023'', funded by Spain's Ministerio de Ciencia, Innovación y Universidades.

\bibliographystyle{apalike-custom}
\bibliography{GMB-arXiv}

\begin{thebibliography}{}

\bibitem[Abrahams and Thomas, 1981]{Abrahams_1981_some}
Abrahams, J. and Thomas, J. (1981).
\newblock Some comments on conditionally {M}arkov and reciprocal {G}aussian processes (corresp.).
\newblock {\em IEEE Transactions on Information Theory}, 27(4):523--525.
\newblock \href {http://dx.doi.org/10.1109/TIT.1981.1056361} {\path{doi:10.1109/TIT.1981.1056361}}.

\bibitem[Andersson, 2012]{Andersson_2012_card}
Andersson, P. (2012).
\newblock Card counting in continuous time.
\newblock {\em Journal of Applied Probability}, 49(1):184--198.
\newblock \href {http://dx.doi.org/10.1239/jap/1331216841} {\path{doi:10.1239/jap/1331216841}}.

\bibitem[Angoshtari and Leung, 2019]{Angoshtari_2019_optimal}
Angoshtari, B. and Leung, T. (2019).
\newblock Optimal dynamic basis trading.
\newblock {\em Annals of Finance}, 15(3):307--335.
\newblock \href {http://dx.doi.org/10.1007/s10436-019-00348-x} {\path{doi:10.1007/s10436-019-00348-x}}.

\bibitem[Back, 1992]{Back_1992_insider}
Back, K. (1992).
\newblock Insider trading in continuous time.
\newblock {\em The Review of Financial Studies}, 5(3):387--409.
\newblock \href {http://dx.doi.org/10.1093/rfs/5.3.387} {\path{doi:10.1093/rfs/5.3.387}}.

\bibitem[Barczy and Kern, 2011]{Barczy_2011_general}
Barczy, M. and Kern, P. (2011).
\newblock General alpha-{W}iener bridges.
\newblock {\em Communications on Stochastic Analysis}, 5(3):585--608.
\newblock \href {http://dx.doi.org/10.31390/cosa.5.3.08} {\path{doi:10.31390/cosa.5.3.08}}.

\bibitem[Barczy and Kern, 2013a]{Barczy_2013_representation}
Barczy, M. and Kern, P. (2013a).
\newblock Representations of multidimensional linear process bridges.
\newblock {\em Random Operators and Stochastic Equations}, 21(2):159--189.
\newblock \href {http://dx.doi.org/10.1515/rose-2013-0009} {\path{doi:10.1515/rose-2013-0009}}.

\bibitem[Barczy and Kern, 2013b]{Barczy_2013_sample}
Barczy, M. and Kern, P. (2013b).
\newblock Sample path deviations of the {W}iener and the {O}rnstein-{U}hlenbeck process from its bridges.
\newblock {\em Brazilian Journal of Probability and Statistics}, 27(4):437--466.
\newblock \href {http://dx.doi.org/10.1214/11-BJPS175} {\path{doi:10.1214/11-BJPS175}}.

\bibitem[Borisov, 1983]{Borisov_1983_criterion}
Borisov, I.~S. (1983).
\newblock On a criterion for {G}aussian random processes to be {M}arkovian.
\newblock {\em Theory of Probability \& Its Applications}, 27(4):863--865.
\newblock \href {http://dx.doi.org/10.1137/1127097} {\path{doi:10.1137/1127097}}.

\bibitem[Boyce, 1970]{Boyce_1970_stopping}
Boyce, W.~M. (1970).
\newblock Stopping rules for selling bonds.
\newblock {\em The Bell Journal of Economics and Management Science}, 1(1):27--53.
\newblock \href {http://dx.doi.org/10.2307/3003021} {\path{doi:10.2307/3003021}}.

\bibitem[Brennan and Schwartz, 1990]{Brennan_1990_arbitrage}
Brennan, M.~J. and Schwartz, E.~S. (1990).
\newblock Arbitrage in stock index futures.
\newblock {\em The Journal of Business}, 63(1):S7--S31.
\newblock \href {http://dx.doi.org/10.1086/296491} {\path{doi:10.1086/296491}}.

\bibitem[Buonocore et~al., 2013]{Buonocore_2013_some}
Buonocore, A., Caputo, L., Nobile, A.~G., and Pirozzi, E. (2013).
\newblock On some time-non-homogeneous linear diffusion processes and related bridges.
\newblock {\em Scientiae Mathematicae Japonicae}, 76(1):55--77.
\newblock \href {http://dx.doi.org/10.32219/isms.76.1_55} {\path{doi:10.32219/isms.76.1_55}}.

\bibitem[Campi and Çetin, 2007]{Campi_2007_insider}
Campi, L. and Çetin, U. (2007).
\newblock Insider trading in an equilibrium model with default: a passage from reduced-form to structural modelling.
\newblock {\em Finance and Stochastics}, 11(4):591--602.
\newblock \href {http://dx.doi.org/10.1007/s00780-007-0038-4} {\path{doi:10.1007/s00780-007-0038-4}}.

\bibitem[Campi et~al., 2011]{Campi_2011_dynamic}
Campi, L., Çetin, U., and Danilova, A. (2011).
\newblock Dynamic {Markov} bridges motivated by models of insider trading.
\newblock {\em Stochastic Processes and their Applications}, 121(3):534--567.
\newblock \href {http://dx.doi.org/10.1016/j.spa.2010.11.004} {\path{doi:10.1016/j.spa.2010.11.004}}.

\bibitem[Campi et~al., 2013]{Campi_2013_equilibrium}
Campi, L., Çetin, U., and Danilova, A. (2013).
\newblock Equilibrium model with default and dynamic insider information.
\newblock {\em Finance and Stochastics}, 17(3):565--585.
\newblock \href {http://dx.doi.org/10.1007/s00780-012-0196-x} {\path{doi:10.1007/s00780-012-0196-x}}.

\bibitem[Cartea et~al., 2016]{Cartea_2016_algorithmic}
Cartea, {\'A}., Jaimungal, S., and Kinzebulatov, D. (2016).
\newblock Algorithmic trading with learning.
\newblock {\em International Journal of Theoretical and Applied Finance}, 19(04):1650028.
\newblock \href {http://dx.doi.org/10.1142/S021902491650028X} {\path{doi:10.1142/S021902491650028X}}.

\bibitem[Cetin and Xing, 2013]{Cetin_2013_point}
Cetin, U. and Xing, H. (2013).
\newblock Point process bridges and weak convergence of insider trading models.
\newblock {\em Electronic Journal of Probability}, 18(26):1--24.
\newblock \href {http://dx.doi.org/10.1214/EJP.v18-2039} {\path{doi:10.1214/EJP.v18-2039}}.

\bibitem[Chaumont and Bravo, 2011]{Chaumont_2011_Markovian}
Chaumont, L. and Bravo, G.~U. (2011).
\newblock Markovian bridges: Weak continuity and pathwise constructions.
\newblock {\em The Annals of Probability}, 39(2):609--647.
\newblock \href {http://dx.doi.org/10.1214/10-AOP562} {\path{doi:10.1214/10-AOP562}}.

\bibitem[Chen et~al., 2015]{Chen_2015_optimal}
Chen, R.~W., Grigorescu, I., and Kang, M. (2015).
\newblock Optimal stopping for {Shepp}'s urn with risk aversion.
\newblock {\em Stochastics. An International Journal of Probability and Stochastic Processes}, 87(4):702--722.
\newblock \href {http://dx.doi.org/10.1080/17442508.2014.995660} {\path{doi:10.1080/17442508.2014.995660}}.

\bibitem[Chen et~al., 2021]{Chen_2021_constrained}
Chen, X., Leung, T., and Zhou, Y. (2021).
\newblock Constrained dynamic futures portfolios with stochastic basis.
\newblock {\em Annals of Finance}, pp. 1--33.
\newblock \href {http://dx.doi.org/10.1007/s10436-021-00398-0} {\path{doi:10.1007/s10436-021-00398-0}}.

\bibitem[Chen and Georgiou, 2016]{Chen_2016_stochastic}
Chen, Y. and Georgiou, T. (2016).
\newblock Stochastic bridges of linear systems.
\newblock {\em IEEE Transactions on Automatic Control}, 61(2):526--531.
\newblock \href {http://dx.doi.org/10.1109/TAC.2015.2440567} {\path{doi:10.1109/TAC.2015.2440567}}.

\bibitem[D'Auria et~al., 2021]{DAuria-2021-optimal}
D'Auria, B., Garc{\'\i}a-Portugu{\'e}s, E., and Guada, A. (2021).
\newblock Optimal stopping of an {O}rnstein--{U}hlenbeck bridge.
\newblock {\em arXiv:2110.13056}.
\newblock \href {http://dx.doi.org/10.48550/arXiv.2110.13056} {\path{doi:10.48550/arXiv.2110.13056}}.

\bibitem[De~Angelis, 2015]{DeAngelis_2015_note}
De~Angelis, T. (2015).
\newblock A note on the continuity of free-boundaries in finite-horizon optimal stopping problems for one-dimensional diffusions.
\newblock {\em SIAM Journal on Control and Optimization}, 53(1):167--184.
\newblock \href {http://dx.doi.org/10.1137/130920472} {\path{doi:10.1137/130920472}}.

\bibitem[De~Angelis and Milazzo, 2020]{DeAngelis_2020_optimal}
De~Angelis, T. and Milazzo, A. (2020).
\newblock Optimal stopping for the exponential of a {B}rownian bridge.
\newblock {\em Journal of Applied Probability}, 57(1):361--384.
\newblock \href {http://dx.doi.org/10.1017/jpr.2019.98} {\path{doi:10.1017/jpr.2019.98}}.

\bibitem[De~Angelis and Peskir, 2020]{DeAngelis_2020_global}
De~Angelis, T. and Peskir, G. (2020).
\newblock Global {$C^{1}$} regularity of the value function in optimal stopping problems.
\newblock {\em The Annals of Applied Probability}, 30(3):1007--1031.
\newblock \href {http://dx.doi.org/10.1214/19-aap1517} {\path{doi:10.1214/19-aap1517}}.

\bibitem[De~Angelis and Stabile, 2019]{DeAngelis_2019_Lipschitz}
De~Angelis, T. and Stabile, G. (2019).
\newblock On {L}ipschitz continuous optimal stopping boundaries.
\newblock {\em SIAM Journal on Control and Optimization}, 57(1):402--436.
\newblock \href {http://dx.doi.org/10.1137/17m1113709} {\path{doi:10.1137/17m1113709}}.

\bibitem[Detemple and Kitapbayev, 2020]{Detemple_2020_value}
Detemple, J. and Kitapbayev, Y. (2020).
\newblock The value of green energy under regulation uncertainty.
\newblock {\em Energy Economics}, 89:104807.
\newblock \href {http://dx.doi.org/10.1016/j.eneco.2020.104807} {\path{doi:10.1016/j.eneco.2020.104807}}.

\bibitem[Dochviri, 1995]{Dochviri_1995_optimal}
Dochviri, B. (1995).
\newblock On optimal stopping of inhomogeneous standard {Markov} processes.
\newblock {\em Georgian Mathematical Journal}, 2(4):335--346.
\newblock \href {http://dx.doi.org/10.1007/BF02255984} {\path{doi:10.1007/BF02255984}}.

\bibitem[{Dynkin}, 1963]{Dynkin_1963_optimum}
{Dynkin}, E.~B. (1963).
\newblock The optimum choice of the instant for stopping a {Markov} process.
\newblock {\em Soviet Mathematics. Doklady}, 150(2):627--629.

\bibitem[D’Auria and Ferriero, 2020]{DAuria_2020_class}
D’Auria, B. and Ferriero, A. (2020).
\newblock A class of {Itô} diffusions with known terminal value and specified optimal barrier.
\newblock {\em Mathematics}, 8(1):123.
\newblock \href {http://dx.doi.org/10.3390/math8010123} {\path{doi:10.3390/math8010123}}.

\bibitem[D’Auria et~al., 2020]{DAuria_2020_discounted}
D’Auria, B., Garc\'ia-Portugu\'es, E., and Guada, A. (2020).
\newblock Discounted optimal stopping of a {B}rownian bridge, with application to {A}merican options under pinning.
\newblock {\em Mathematics}, 8(7):1159.
\newblock \href {http://dx.doi.org/10.3390/math8071159} {\path{doi:10.3390/math8071159}}.

\bibitem[Ekström and Vaicenavicius, 2020]{Ekstrom_2020_optimal}
Ekström, E. and Vaicenavicius, J. (2020).
\newblock Optimal stopping of a {Brownian} bridge with an unknown pinning point.
\newblock {\em Stochastic Processes and their Applications}, 130(2):806--823.
\newblock \href {http://dx.doi.org/10.1016/j.spa.2019.03.018} {\path{doi:10.1016/j.spa.2019.03.018}}.

\bibitem[Ekström and Wanntorp, 2009]{Ekstrom_2009_optimal}
Ekström, E. and Wanntorp, H. (2009).
\newblock Optimal stopping of a {Brownian} bridge.
\newblock {\em Journal of Applied Probability}, 46(1):170--180.
\newblock \href {http://dx.doi.org/10.1239/jap/1238592123} {\path{doi:10.1239/jap/1238592123}}.

\bibitem[Erickson and Steck, 2022]{Erickson_2020_probability}
Erickson, W.~W. and Steck, D.~A. (2022).
\newblock The anatomy of an extreme event: What can we infer about the history of a heavy-tailed random walk?
\newblock {\em arXiv:2002.03849}.
\newblock \href {http://dx.doi.org/10.48550/arXiv.2002.03849} {\path{doi:10.48550/arXiv.2002.03849}}.

\bibitem[Ernst and Shepp, 2015]{Ernst_2015_revisiting}
Ernst, P.~A. and Shepp, L.~A. (2015).
\newblock Revisiting a theorem of {L}. {A}. {Shepp} on optimal stopping.
\newblock {\em Communications on Stochastic Analysis}, 9(3):419--423.
\newblock \href {http://dx.doi.org/10.31390/cosa.9.3.08} {\path{doi:10.31390/cosa.9.3.08}}.

\bibitem[Fitzsimmons et~al., 1993]{Fitzsimmons_1993_Markovian}
Fitzsimmons, P., Pitman, J., and Yor, M. (1993).
\newblock Markovian bridges: Construction, palm interpretation, and splicing.
\newblock In Çinlar, E., Chung, K.~L., Sharpe, M.~J., Bass, R.~F., and Burdzy, K. (Eds.), {\em Seminar on Stochastic Processes, 1992}, volume~33 of {\em Progress in Probability}, pp. 101--134. Birkh\"auser, Boston.
\newblock \href {http://dx.doi.org/10.1007/978-1-4612-0339-1_5} {\path{doi:10.1007/978-1-4612-0339-1_5}}.

\bibitem[Friedman, 1964]{Friedman_1964_partial}
Friedman, A. (1964).
\newblock {\em Partial Differential Equations of Parabolic Type}.
\newblock Prentice-Hall, Englewood Cliffs.

\bibitem[Friedman, 1975a]{Friedman_1975_parabolic}
Friedman, A. (1975a).
\newblock Parabolic variational inequalities in one space dimension and smoothness of the free boundary.
\newblock {\em Journal of Functional Analysis}, 18(2):151--176.
\newblock \href {http://dx.doi.org/10.1016/0022-1236(75)90022-1} {\path{doi:10.1016/0022-1236(75)90022-1}}.

\bibitem[Friedman, 1975b]{Friedman_1975_stopping}
Friedman, A. (1975b).
\newblock Stopping {Time} {Problems} and the {Shape} of the {Domain} of {Continuation}.
\newblock In {\em Control {Theory}, {Numerical} {Methods} and {Computer} {Systems} {Modelling}}, Lecture {Notes} in {Economics} and {Mathematical} {Systems}, pp. 559--566, Berlin, Heidelberg. Springer.
\newblock \href {http://dx.doi.org/10.1007/978-3-642-46317-4_39} {\path{doi:10.1007/978-3-642-46317-4_39}}.

\bibitem[Föllmer, 1972]{Follmer_1972_optimal}
Föllmer, H. (1972).
\newblock Optimal stopping of constrained {Brownian} motion.
\newblock {\em Journal of Applied Probability}, 9(3):557--571.
\newblock \href {http://dx.doi.org/10.2307/3212325} {\path{doi:10.2307/3212325}}.

\bibitem[Gasbarra et~al., 2007]{Gasbarra_2007_Gaussian}
Gasbarra, D., Sottinen, T., and Valkeila, E. (2007).
\newblock Gaussian bridges.
\newblock In Benth, F.~E., Di~Nunno, G., Lindstr\o{}m, T., \O{}ksendal, B., and Zhang, T. (Eds.), {\em Stochastic Analysis and Applications}, Abel Symposia, pp. 361--382, Berlin. Springer.
\newblock \href {http://dx.doi.org/10.1007/978-3-540-70847-6_15} {\path{doi:10.1007/978-3-540-70847-6_15}}.

\bibitem[Glover, 2020]{Glover_2020_optimally}
Glover, K. (2020).
\newblock Optimally stopping a {Brownian} bridge with an unknown pinning time: a {Bayesian} approach.
\newblock {\em Stochastic Processes and their Applications}, 150:919--937.
\newblock \href {http://dx.doi.org/10.1016/j.spa.2020.03.007} {\path{doi:10.1016/j.spa.2020.03.007}}.

\bibitem[Golez and Jackwerth, 2012]{Golez_2012_pinning}
Golez, B. and Jackwerth, J.~C. (2012).
\newblock Pinning in the {S}\&{P} 500 futures.
\newblock {\em Journal of Financial Economics}, 106(3):566--585.
\newblock \href {http://dx.doi.org/10.1016/j.jfineco.2012.06.010} {\path{doi:10.1016/j.jfineco.2012.06.010}}.

\bibitem[Hildebrandt and Rœlly, 2020]{Hildebrandt_2020_pinned}
Hildebrandt, F. and Rœlly, S. (2020).
\newblock Pinned diffusions and {M}arkov bridges.
\newblock {\em Journal of Theoretical Probability}, 33(2):906--917.
\newblock \href {http://dx.doi.org/10.1007/s10959-019-00954-5} {\path{doi:10.1007/s10959-019-00954-5}}.

\bibitem[Hilliard and Hilliard, 2015]{Hilliard_2015_pricing}
Hilliard, J.~E. and Hilliard, J. (2015).
\newblock Pricing {American} options when there is short--lived arbitrage.
\newblock {\em International Journal of Financial Markets and Derivatives}, 4(1):43--53.
\newblock \href {http://dx.doi.org/10.1504/IJFMD.2015.066444} {\path{doi:10.1504/IJFMD.2015.066444}}.

\bibitem[Horne et~al., 2007]{Horne_2007_analyzing}
Horne, J.~S., Garton, E.~O., Krone, S.~M., and Lewis, J.~S. (2007).
\newblock Analyzing animal movements using {B}rownian bridges.
\newblock {\em Ecology}, 88(9):2354--2363.
\newblock \href {http://dx.doi.org/10.1890/06-0957.1} {\path{doi:10.1890/06-0957.1}}.

\bibitem[Hoyle et~al., 2011]{Hoyle_2011_Levy}
Hoyle, E., Hughston, L.~P., and Macrina, A. (2011).
\newblock L\'evy random bridges and the modelling of financial information.
\newblock {\em Stochastic Processes and Their Applications}, 121(4):856--884.
\newblock \href {http://dx.doi.org/10.1016/j.spa.2010.12.003} {\path{doi:10.1016/j.spa.2010.12.003}}.

\bibitem[Jacka and Lynn, 1992]{Jacka_1992_finite-horizon}
Jacka, S. and Lynn, R. (1992).
\newblock Finite-horizon optimal stopping, obstacle problems and the shape of the continuation region.
\newblock {\em Stochastics and Stochastics Reports}, 39(1):25--42.
\newblock \href {http://dx.doi.org/10.1080/17442509208833761} {\path{doi:10.1080/17442509208833761}}.

\bibitem[Kranstauber, 2019]{Kranstauber_2019_modelling}
Kranstauber, B. (2019).
\newblock Modelling animal movement as {Brownian} bridges with covariates.
\newblock {\em Movement Ecology}, 7(1):22.
\newblock \href {http://dx.doi.org/10.1186/s40462-019-0167-3} {\path{doi:10.1186/s40462-019-0167-3}}.

\bibitem[Krishnan and Nelken, 2001]{Nelken_2001_effect}
Krishnan, H. and Nelken, I. (2001).
\newblock The effect of stock pinning upon option prices.
\newblock {\em Risk}, December:17--20.

\bibitem[Krumm, 2021]{Krumm_2021_Brownian}
Krumm, J. (2021).
\newblock Brownian bridge interpolation for human mobility?
\newblock In {\em Proceedings of the 29th {International} {Conference} on {Advances} in {Geographic} {Information} {Systems}}, {SIGSPATIAL} '21, pp. 175--183, New York, USA. Association for Computing Machinery.
\newblock \href {http://dx.doi.org/10.1145/3474717.3483942} {\path{doi:10.1145/3474717.3483942}}.

\bibitem[Krylov and Aries, 1980]{Krylov_1980_controlled}
Krylov, N.~V. and Aries, A.~B. (1980).
\newblock {\em Controlled Diffusion Processes}.
\newblock Stochastic Modelling and Applied Probability. Springer, New York.

\bibitem[Kyle, 1985]{Kyle_1985_continuous}
Kyle, A.~S. (1985).
\newblock Continuous auctions and insider trading.
\newblock {\em Econometrica}, 53(6):1315--1335.
\newblock \href {http://dx.doi.org/10.2307/1913210} {\path{doi:10.2307/1913210}}.

\bibitem[Leung et~al., 2018]{Leung_2018_optimal}
Leung, T., Li, J., and Li, X. (2018).
\newblock Optimal timing to trade along a randomized {Brownian} bridge.
\newblock {\em International Journal of Financial Studies}, 6(3).
\newblock \href {http://dx.doi.org/10.3390/ijfs6030075} {\path{doi:10.3390/ijfs6030075}}.

\bibitem[Liu and Longstaff, 2004]{Liu_2004_losing}
Liu, J. and Longstaff, F.~A. (2004).
\newblock Losing money on arbitrage: optimal dynamic portfolio choice in markets with arbitrage opportunities.
\newblock {\em The Review of Financial Studies}, 17(3):611--641.
\newblock \href {http://dx.doi.org/10.2139/ssrn.246835} {\path{doi:10.2139/ssrn.246835}}.

\bibitem[Mehr and McFadden, 1965]{Mehr_1965_certain}
Mehr, C.~B. and McFadden, J.~A. (1965).
\newblock Certain properties of {G}aussian processes and their first-passage times.
\newblock {\em Journal of the Royal Statistical Society, Series B (Methodological)}, 27(3):505--522.
\newblock \href {http://dx.doi.org/10.1111/j.2517-6161.1965.tb00611.x} {\path{doi:10.1111/j.2517-6161.1965.tb00611.x}}.

\bibitem[Ni et~al., 2005]{Ni_2005_stock}
Ni, S.~X., Pearson, N.~D., and Poteshman, A.~M. (2005).
\newblock Stock price clustering on option expiration dates.
\newblock {\em Journal of Financial Economics}, 78(1):49--87.
\newblock \href {http://dx.doi.org/10.1016/j.jfineco.2004.08.005} {\path{doi:10.1016/j.jfineco.2004.08.005}}.

\bibitem[Ni et~al., 2021]{Ni_2021_does}
Ni, S.~X., Pearson, N.~D., Poteshman, A.~M., and White, J. (2021).
\newblock Does option trading have a pervasive impact on underlying stock prices?
\newblock {\em The Review of Financial Studies}, 34(4):1952--1986.
\newblock \href {http://dx.doi.org/10.1093/rfs/hhaa082} {\path{doi:10.1093/rfs/hhaa082}}.

\bibitem[Oshima, 2006]{Oshima_2006_optimal}
Oshima, Y. (2006).
\newblock On an optimal stopping problem of time inhomogeneous diffusion processes.
\newblock {\em SIAM Journal on Control and Optimization}, 45(2):565--579.
\newblock \href {http://dx.doi.org/10.1137/040609549} {\path{doi:10.1137/040609549}}.

\bibitem[Pedersen and Peskir, 2002]{Pedersen_2002_onnonlinear}
Pedersen, J.~L. and Peskir, G. (2002).
\newblock On nonlinear integral equations arising in problems of optimal stopping.
\newblock In Baki{\'c}, D., Pand\v{z}i{\'c}, P., and Peskir, G. (Eds.), {\em Functional analysis VII: Proceedings of the Postgraduate School and Conference held in Dubrovnik, September 17-26, 2001}, volume~46 of {\em Various publications series}, pp. 159--175. University of Aarhus, Department of Mathematical Sciences, Aarhus.

\bibitem[Peng and Zhu, 2006]{Peng_2006_necessary}
Peng, S. and Zhu, X. (2006).
\newblock Necessary and sufficient condition for comparison theorem of 1-dimensional stochastic differential equations.
\newblock {\em Stochastic Processes and their Applications}, 116(3):370--380.
\newblock \href {http://dx.doi.org/10.1016/j.spa.2005.08.004} {\path{doi:10.1016/j.spa.2005.08.004}}.

\bibitem[Peskir, 2005a]{Peskir_2005_change}
Peskir, G. (2005a).
\newblock A change-of-variable formula with local time on curves.
\newblock {\em Journal of Theoretical Probability}, 18(3):499--535.
\newblock \href {http://dx.doi.org/10.1007/s10959-005-3517-6} {\path{doi:10.1007/s10959-005-3517-6}}.

\bibitem[Peskir, 2005b]{Peskir_2005_American}
Peskir, G. (2005b).
\newblock On the {A}merican option problem.
\newblock {\em Mathematical Finance}, 15(1):169--181.
\newblock \href {http://dx.doi.org/10.1111/j.0960-1627.2005.00214.x} {\path{doi:10.1111/j.0960-1627.2005.00214.x}}.

\bibitem[Peskir, 2019]{Peskir_2019_continuity}
Peskir, G. (2019).
\newblock Continuity of the optimal stopping boundary for two-dimensional diffusions.
\newblock {\em The Annals of Applied Probability}, 29(1):505--530.
\newblock \href {http://dx.doi.org/10.1214/18-aap1426} {\path{doi:10.1214/18-aap1426}}.

\bibitem[Peskir and Shiryaev, 2006]{Peskir_2006_optimal}
Peskir, G. and Shiryaev, A. (2006).
\newblock {\em Optimal Stopping and Free-Boundary Problems}.
\newblock Lectures in Mathematics. ETH Z\"urich. Birkh\"auser, Basel.
\newblock \href {http://dx.doi.org/10.1007/978-3-7643-7390-0} {\path{doi:10.1007/978-3-7643-7390-0}}.

\bibitem[Pitman and Yor, 1982]{Pitman_1982_decomposition}
Pitman, J. and Yor, M. (1982).
\newblock A decomposition of {B}essel bridges.
\newblock {\em Zeitschrift für Wahrscheinlichkeitstheorie und Verwandte Gebiete}, 59(4):425--457.
\newblock \href {http://dx.doi.org/10.1007/BF00532802} {\path{doi:10.1007/BF00532802}}.

\bibitem[Rosén, 1965]{Rosen_1965_limit}
Rosén, B. (1965).
\newblock Limit theorems for sampling from finite populations.
\newblock {\em Arkiv för Matematik}, 5(5):383--424.
\newblock \href {http://dx.doi.org/10.1007/BF02591138} {\path{doi:10.1007/BF02591138}}.

\bibitem[Salminen, 1984]{Salminen_1984_Brownian}
Salminen, P. (1984).
\newblock Brownian excursions revisited.
\newblock In Çinlar, E., Chung, K.~L., and Getoor, R.~K. (Eds.), {\em Seminar on Stochastic Processes, 1983}, volume~7 of {\em Progress in Probability and Statistics}, pp. 161--187. Birkh\"auser, Boston.
\newblock \href {http://dx.doi.org/10.1007/978-1-4684-9169-2_11} {\path{doi:10.1007/978-1-4684-9169-2_11}}.

\bibitem[Shepp, 1969]{Shepp_1969_explicit}
Shepp, L.~A. (1969).
\newblock Explicit solutions to some problems of optimal stopping.
\newblock {\em Annals of Mathematical Statistics}, 40(3):993--1010.

\bibitem[Shiryaev, 2008]{Shiryaev_2008_optimal}
Shiryaev, A. (2008).
\newblock {\em Optimal {Stopping} {Rules}}.
\newblock Stochastic {Modelling} and {Applied} {Probability}. Springer-Verlag, Berlin Heidelberg.
\newblock \href {http://dx.doi.org/10.1007/978-3-540-74011-7} {\path{doi:10.1007/978-3-540-74011-7}}.

\bibitem[Sottinen and Yazigi, 2014]{Sottinen_2014_generalized}
Sottinen, T. and Yazigi, A. (2014).
\newblock Generalized {Gaussian} bridges.
\newblock {\em Stochastic Processes and their Applications}, 124(9):3084--3105.
\newblock \href {http://dx.doi.org/10.1016/j.spa.2014.04.002} {\path{doi:10.1016/j.spa.2014.04.002}}.

\bibitem[Taylor, 1968]{Taylor_1968_optimal}
Taylor, H.~M. (1968).
\newblock Optimal stopping in a {Markov} process.
\newblock {\em Annals of Mathematical Statistics}, 39(4):1333--1344.
\newblock \href {http://dx.doi.org/10.1214/aoms/1177698259} {\path{doi:10.1214/aoms/1177698259}}.

\bibitem[Venek et~al., 2016]{Venek_2016_evaluating}
Venek, V., Brunauer, R., and Schneider, C. (2016).
\newblock Evaluating the {Brownian} bridge movement model to determine regularities of people’s movements.
\newblock {\em Journal for Geographic Information Science}, 4:20--35.
\newblock \href {http://dx.doi.org/10.1553/giscience2016_02_s20} {\path{doi:10.1553/giscience2016_02_s20}}.

\bibitem[Wald, 1947]{Wald_1947_sequential}
Wald, A. (1947).
\newblock {\em Sequential Analysis}.
\newblock Wiley Series in Probability and Mathematical Statistics. John Wiley \& Sons, New York.

\bibitem[Williams and Rasmussen, 2006]{Ramussen_2006_Gaussian}
Williams, C.~K. and Rasmussen, C.~E. (2006).
\newblock {\em Gaussian Processes for Machine Learning}.
\newblock MIT press, Cambridge.
\newblock \href {http://dx.doi.org/10.7551/mitpress/3206.001.0001} {\path{doi:10.7551/mitpress/3206.001.0001}}.

\bibitem[Yang, 2014]{Yang_2014_refined}
Yang, Y. (2014).
\newblock Refined solutions of time inhomogeneous optimal stopping problem and zero-sum game via {Dirichlet} form.
\newblock {\em Probability and Mathematical Statistics}, 34(2):253--271.

\bibitem[Çetin and Danilova, 2016]{Cetin_2016_Markov}
Çetin, U. and Danilova, A. (2016).
\newblock Markov bridges: {SDE} representation.
\newblock {\em Stochastic Processes and Their Applications}, 126(3):651--679.
\newblock \href {http://dx.doi.org/10.1016/j.spa.2015.09.015} {\path{doi:10.1016/j.spa.2015.09.015}}.

\end{thebibliography}

\end{document}